\documentclass[11pt,reqno]{amsart}

\usepackage[text={160mm,240mm},centering]{geometry}            
\usepackage{amssymb,amsmath,setspace,geometry,indentfirst,changepage,inputenc,amsthm}
\usepackage{mathrsfs,amsfonts,savesym,graphicx,bm,color}

\usepackage{graphicx}
\usepackage{amssymb}
\usepackage{epstopdf}
\usepackage{dsfont}

\usepackage{lscape}

\usepackage{tikz}
\usetikzlibrary{positioning}
\usepackage{xcolor}

\usepackage{float}

\usepackage[all,2cell]{xy}    
\UseAllTwocells               

\usepackage{amsmath,amsthm}

\usepackage[mathscr]{eucal}
\usepackage[all]{xy}
\usepackage{mathrsfs}
\usepackage{hyperref}
\usepackage{color}

\newtheorem{theorem}{Theorem}[section]
\newtheorem{lemma}[theorem]{Lemma}

\newtheorem{thmletter}{Theorem}

\newtheorem{conjecture}[theorem]{Conjecture}
\newtheorem{corollary}[theorem]{Corollary}
\newtheorem{proposition}[theorem]{Proposition}

\newtheorem{definition-lemma}[theorem]{Definition-Lemma}
\newtheorem{definition-theorem}[theorem]{Definition-Theorem}

\theoremstyle{definition}
\newtheorem{example}[theorem]{Example}
\newtheorem{definition}[theorem]{Definition}

\newtheorem{remark}[theorem]{Remark}

\allowdisplaybreaks

\title[]{On jet schemes of determinantal varieties}

\author{Yifan Chen}
\address{Department of Mathematical Sciences,
	Tsinghua University,
	Beijing, 100084, P. R. China.}
\email{c-yf20@tsinghua.org.cn}
\author{Huaiqing Zuo}
\address{Department of Mathematical Sciences,
	Tsinghua University,
	Beijing, 100084, P. R. China.}
\email{hqzuo@mail.tsinghua.edu.cn}

\begin{document}
	
	\maketitle
	
	\parskip=.3em
	
	\setcounter{tocdepth}{1}

\begin{abstract}
Determinantal varieties are important objects of study in algebraic geometry. In this paper, we will investigate them using the jet scheme approach. We have found a new connection for the Hilbert series between a determinantal variety and its jet schemes. We denote the $k$-th order jet scheme of the determinantal variety defined by $r$-minors in an $m \times n$ matrix as $\mathscr{L}^{m,n}_{r,k}$. For the special case where $m$, $n$, and $r$ are equal, and $m$ and $r$ are 3 while $k$ is 1, we establish a correspondence between the defining ideals of $\mathscr{L}^{m,n}_{r,k}$ and abstract simplicial complexes, proving their shellability and obtaining the Hilbert series of $\mathscr{L}^{m,n}_{r,k}$ accordingly. Moreover, for general $\mathscr{L}^{m,n}_{r,k}$, \cite{12} provides its irreducible decomposition. We further provide a specific polynomial family defining its irreducible components.
	

\end{abstract}

\section{introduction}
In recent years, the theory of jet schemes and arc spaces has been widely applied. They were initially introduced by Nash in \cite{8} in 1968. Subsequently, Denef and Loeser utilized jet schemes and arc spaces for motivic integration in \cite{16}, a concept introduced by Maxim Kontsevich in \cite{9}. These tools were also employed by Mustaţă in \cite{10} to study the singularities of pairs $(X,Y)$, where $X$ is a smooth algebraic variety and $Y \subset X$ is a closed subscheme. Ein and Mustaţă also utilized arc spaces and jet schemes to study other invariants of singularities, some of which are summarized in \cite{11}.

Let $K$ be an algebraically closed field, $\mathbf{A}_{K}^{n}$ be an affine space of dimension $n$. For a fixed $m \times n$ matrix (without loss of generality, we can assume $m \le n$), let $\mathbf{A}_{K}^{mn}$ be its affine space, and define an algebraic variety $\mathscr{L}^{m,n}_{r}$ in this space determined by all $r \times r$ minors of the matrix. This algebraic variety $\mathscr{L}^{m,n}_{r}$ is called the determinantal variety, which is an important object in the study of algebraic geometry. In this article, we denote the variables of an $m \times n$ matrix as $x_{i,j}(1 \le i \le m, 1 \le j \le n)$ and the ideal generated by the $r \times r$ minors as $I^{m,n}_{r}$. Then, the coordinate ring of $\mathscr{L}^{m,n}_{r}$ is $K[x_{i,j}|1 \le i \le m,1 \le j \le n]/I^{m,n}_{r}$.

In the past decade, there has been considerable research on the jet schemes of determinantal varieties. When computing the jet scheme of $\mathscr{L}^{m,n}_{r,k}$, we assume $x_{i,j}=x_{i,j}^{(0)}+x_{i,j}^{(1)}t+...+x_{i,j}^{(k)}t^{k}$. The coordinate ring of the resulting jet scheme is given by $K[x_{i,j}^{(l)}|1 \le i \le m, 1 \le j \le n, 0 \le l \le k]/I^{m,n}_{r,k}$, where $I^{m,n}_{r,k}$ is the jet ideal of $I^{m,n}_{r}$ (see Section 2.1 for the computation of jet schemes). In \cite{12}, Roi Docampo investigated the irreducible components of the jet schemes of determinantal varieties, computed their number, and provided a characterization. The central result can be stated as the following theorem:

\begin{theorem}
(\cite{12}, Proposition $4.9$, Proposition $4.11$, Corollary $4.13$)Every irreducible components of $\mathscr{L}^{m,n}_{r,p-1}$ has a one to one correspondence with a pre-partition $\lambda=(d \ d \ ... \ d \ e \ 0 \ ... \ 0)$, satisfying $p=(a+1)d+e$, $0 \le e < d$; either $e=0, 0 \le a \le r-1$, or $e > 0, 0 \le a < r-1$, and there are totally $a+m-r+1$ many $d$ in $\lambda$. By calculating the number of $\lambda$ we can get: when $r=1,m$, $\mathscr{L}^{m,n}_{r,k}$ is irreducible; when $1<r<m$, $\mathscr{L}^{m,n}_{r,k}$ has $k+2-\lceil\frac{k+1}{r}\rceil$ many irreducible components.
\end{theorem}

For general $m,n,r$ and $k$, the computation of each irreducible component of $\mathscr{L}^{m,n}_{r,k}$ is not easy, but there have been quite a few results in this regard, especially in the case of small parameters. The simplest non-trivial case is when $r=2$ and $k=1$, which corresponds to the first jet scheme of determinantal varieties defined by $2 \times 2$ minors. According to the above theorem, in this case, $\mathscr{L}^{m,n}_{2,1}$ has two irreducible components. Let us denote the ideals corresponding to these two irreducible components as $I_{1}$ and $I_{2}$. In \cite{2}, the authors give explicit expressions for them and prove that $I^{m,n}_{2,1}=I_{1} \cap I_{2}$, which is precisely the primary decomposition of $I^{m,n}_{2,1}$. Among these two components, one corresponds to the closure of the smooth part of the base space $\mathscr{L}^{m,n}_{2}$ under the inverse image of the projection map $\mathscr{L}^{m,n}_{2,1} \longrightarrow \mathscr{L}^{m,n}_{2}$. This component is called the principal component, and we denote its corresponding ideal as $I_{1}$. In \cite{2}, the authors directly provide a set of Gr\"{o}bner bases for both $I^{m,n}_{2,1}$ and $I_{1}$. Then, in \cite{13}, Boyan Jonov characterizes the Stanley-Reisner abstract simplicial complex corresponding to the initial ideal of $I_{1}$ and proves its shellability(see Definition $3.5$). Finally, in \cite{5}, the authors use the above results to compute the Hilbert series of the principal component of $\mathscr{L}^{m,n}_{2,1}$.

In addition to the above cases, there have been some good results for the situation where $r=m=n$. Firstly, according to Theorem $1.1$, when $r=m$, $\mathscr{L}^{m,n}_{m,k}$ is irreducible. If we add the condition $n=m$, in \cite{14}, the authors proved that under a certain defined monomial order, the polynomials defining $\mathscr{L}^{m,m}_{m,k}$ exactly form a Gr\"{o}bner bases for the ideal they generate.

Building on previous work, this paper delves deeper into the study of more general and larger parameter settings. Firstly, when $r=m=n$, we prove that the Stanley-Reisner abstract simplicial complex(see Definition $3.4$) corresponding to the initial ideal of $I^{m,m}_{m,k}$ is shellable, and we have calculated the Hilbert series of $\mathscr{L}^{m,m}_{m,k}$. The specific description is given in the following theorem.

\begin{thmletter}\label{mtha}
Let $LM(I^{m,m}_{m,k})$ be the leading ideal of $I^{m,m}_{m,k}$, and assume the Stanley-Reisner abstract simplicial complex corresponding to $R=K[x_{i,j}^{(l)}|1 \le i,j \le m, 0 \le l \le k]/LM(I^{m,m}_{m,k})$ is $\Delta$, then $\Delta$ is shellable and the Hilbert series of $R$ is as following($z$ is the variable):
\begin{flalign*}
	\frac{(1-z^{m})^{k+1}}{(1-z)^{m^{2}(k+1)}}.
\end{flalign*}
\end{thmletter}

When $m<n$, the general case becomes more challenging. This paper specifically focuses on the case where $m=3$ and $k=1$. It proves that the abstract simplicial complex corresponding to $LM(I^{3,n}_{3,1})$ is shellable and uses this result to compute the Hilbert series of $\mathscr{L}^{3,n}_{3,1}$, as shown in the following theorem.
\begin{thmletter}\label{mthb}
Suppose the Stanley-Reisner abstract simplicial complex corresponding to $R=K[x_{i,j}^{(l)}|1 \le i \le 3, 1 \le j \le n, 0 \le l \le 1]/LM(I^{3,n}_{3,1})$ is $\Delta_{0}$, then $\Delta_{0}$ is shellable and the Hilbert series of $R$ is($z$ is the variable):
\begin{flalign*}
(\frac{1+(n-2)z+\frac{(n-1)(n-2)}{2}z^{2}}{(1-z)^{2n+2}})^{2}.
\end{flalign*}
Moreover, it is exactly the square of the Hilbert series of $K[x_{i,j}|1 \le i \le 3, 1 \le j \le n]/I^{3,n}_{3}$.
\end{thmletter}

In \cite{5}, authors have calculated the Hilbert series of principal component of $\mathscr{L}^{m,n}_{2,1}$, which turns out to be the square of the Hilbert series of the base space $\mathscr{L}^{m,n}_{2}$. The theorem above also shows this result holds for the case $\mathscr{L}^{3,n}_{3,1}$. This shows some specific connection between the base space and the principal component of its jet scheme. For general case $\mathscr{L}^{m,n}_{r,k}$, in section 8 we will prove that the principal component always exists and give an explicit description of it. There are a lot of examples show that the conclusion still holds for the general case, as the following conjecture:
\begin{conjecture}
Denote $\tilde{I^{m,n}_{r,k}}$ as the ideal corresponding to the principal component of $\mathscr{L}^{m,n}_{r,k}$, then its corresponding abstract simplicial complex is shellable, and the Hilbert series of the principal component of $\mathscr{L}^{m,n}_{r,k}$ is exactly $k+1$-th square of the Hilbert series of $\mathscr{L}^{m,n}_{r}$.
\end{conjecture}   

The above theorems are all under the assumption of $r=m$. When $1<r<m$, it is known from Theorem $1.1$ that $\mathscr{L}^{m,n}_{r,k}$ is reducible. Therefore, we need to better characterize its irreducible components. Building on \cite{12}, this paper provides polynomials whose zero locus is the corresponding irreducible component. First, let's establish some notation. For an irreducible component corresponding to a pre-partition $\lambda=(\lambda_{1},...,\lambda_{m})$ in Theorem $1.1$, we define $i_{s}=\lambda_{m}+...+\lambda_{m-s+1},j_{s}=\lceil\frac{i_{s}}{s}\rceil$. We consider an $m \times n$ matrix with variables $x_{i,j}^{(0)}+x_{i,j}^{(1)}t+...+x_{i,j}^{(p-1)}t^{p-1}$, and denote the polynomial family consisting of $t^{\alpha}$ term of the $s$-th order minors in the variables $x_{i,j}^{(l)}(0 \le l \le p-1)$ as $\tilde{\Theta}_{s,\alpha}$, i.e. for an arbitrary $s$-th order minor whose entry is $x_{i,j}^{(0)}+x_{i,j}^{(1)}t+...+x_{i,j}^{(p-1)}t^{p-1}$, we calculate its determinant as $g_{0}+g_{1}t+...+g_{\alpha}t^{\alpha}+...$, then we take $g_{\alpha}$ into the polynomial family $\tilde{\Theta}_{s,\alpha}$. We have the following theorem.
\begin{thmletter}\label{mthc}
For a fixed $\lambda$, we construct the family of polynomials $\tilde{\Omega_{\lambda}'}$ as the following steps: first we take $\tilde{\Theta}_{1,0},...,\tilde{\Theta}_{1,i_{1}-1} \in \tilde{\Omega_{\lambda}'}$; for every $1 \le s \le m-1$, when $i_{s}+j_{s}<i_{s+1}$, we take $\tilde{\Theta}_{s+1,i_{s}+j_{s}},...,\tilde{\Theta}_{s+1,i_{s+1}-1} \in \tilde{\Omega_{\lambda}'}$. We do this procedure for every $1 \le s \le m-1$ to get the whole family $\tilde{\Omega_{\lambda}'}$. Then this family of polynomials can define the irreducible components of  $\mathscr{L}^{m,n}_{r,p-1}$ which corresponds to $\lambda$.
\end{thmletter}

The paper unfolds in the following order: In the second section, we prepare some foundational knowledge, including the definition, properties, and computations of jet schemes, basic knowledge about Gr\"{o}bner bases, and the Buchberger criterion; In the third section, we explore $\mathscr{L}^{m,m}_{m,k}$ and prove Theorem \ref{mtha}; In sections four through seven, we investigate $\mathscr{L}^{3,n}_{3,1}$ in several steps and prove Theorem \ref{mthb}; Finally, in the eighth section, we delve into the irreducible components of general $\mathscr{L}^{m,n}_{r,k}$ and prove Theorem \ref{mthc}.

\section{preliminaries}

\subsection{Jet Scheme and its Calculation}
Now we will give the definition of jet scheme and introduce the method to calculate it.
\begin{definition}
	Let $X$ be a scheme over $K$, for every $m \in \mathbb{N}$, consider the functor from $K$-schemes to set
	\begin{flalign*}
		Z \mapsto Hom(Z \times_{Spec(K)} Spec (K[t]/(t^{m+1})),X).
	\end{flalign*}
	It was shown in {\color{black}Theorem $2.2$ in \cite{6}} that there is a $K$-scheme represents this functor, which is called the $m$-th jet scheme of $X$, denoted by $X_{m}$, i.e.
	\begin{flalign*}
		Hom(Z,X_{m})=Hom(Z \times_{Spec(K)} Spec (K[t]/(t^{m+1})),X).
	\end{flalign*}
\end{definition}

For $1 \le i \le j$, we can get natural projections between jet schemes $\psi_{i,j}: X_{j} \longrightarrow X_{i}$ induced from the truncation map $K[t]/(t^{j}) \longrightarrow K[t]/(t^{i})$. If we identify $X_{0}$ with $X$, we have natural projection $\pi_{i}: X_{i} \longrightarrow X$. It can be shown that $\{X_{m}\}$ is an inverse system and the inverse limit $X_{\infty}:=\lim \limits_{\leftarrow m}X_{m}$ is called the arc space of $X$, so there are also natural projections $\psi_{i}: X_{\infty} \longrightarrow X_{i}$. {\color{black}Similar to jet scheme, the arc space can also be characterized by functor of points if $X$ is of finite type. More precisely, Proposition $2.13$ in \cite{6} shows that $$Hom(Z,X_{\infty})=Hom(Z \hat{\times}_{Spec(K)} Spec (K[[t]]),X)$$ for any $K$-scheme $Z$, where $Z \hat{\times}_{Spec(K)} Spec (K[[t]])$ is the completion of $Z \times_{Spec(K)} Spec (K[[t]])$ along $Z \times_{Spec(K)} \{0\}$. One can also refer to Corollary $1.2$ in \cite{17} for general statement.} Note that $X_{m}$ is of finite type for each $m \in \mathbb{N}$ if $X$ is of finite type, but $X_{\infty}$ is usually not.

To compute jet schemes explicitly, we introduce the following algebra:
\begin{definition}
Suppose $B$ is a $K$-algebra and $m \in \mathbb{N}$ is an integer, we define a $B$-algebra $HS_{B}^{m} :=B[x^{(i)}]_{x \in B, i=1,2,...,m}/I$, where $I$ is generated by the following set: 	
	
(1)$\{ (x+y)^{(i)}-x^{(i)}-y^{(i)}:x,y \in B,i=1,2,...,m\}$. 
	
(2)$\{a^{(i)}:a \in K,i=1,...,m \}$.
	
(3)$\{(xy)^{(k)}-\displaystyle\sum_{i+j=k}x^{(i)}y^{(j)}:x,y \in B, k=0,...,m\}$,
where we identify $x^{(0)}$ with $x$ for all $x \in B$. 
	
We also define the universal derivation $(d_{0},...,d_{m})$ from $B$ to $HS_{B}^{m}$ by $d_{i}(x)=x^{(i)}$.
\end{definition}

The following proposition in \cite{13} shows that this algebra is the one we need to find.
\begin{proposition}
(Corollary $1.8$ in \cite{13}) 
For any $\phi \in Hom(HS_{B}^{m},R)$, define a morphism from $B$ to $R[t]/(t^{m+1})$ given by $x \mapsto \phi(d_{0}(x))+...+\phi(d_{m}(x))t^{m}$, and this gives a bijection between $Hom(HS_{B}^{m},R)$ and $Hom(B,R[t]/(t^{m+1}))$.
\end{proposition} 

From this proposition we know if we set $X=Spec(B)$, then we have $X_{m}=Spec(HS_{B}^{m})$. When the scheme $X=Spec(B)$, we will denote $HS_{B}^{m}$ by $B_{m}$ for the convenience of notation, and we say $B_{m}$ is the $m$-th jet ring of $B$. In particular, if $R$ is a ring, $I \subset R$ is an ideal and $B=R/I$, we have $B_{m}=R_{m}/I_{m}$, and $I_{m}$ is called $m$-th jet ideal.

We will use determinantal variety as an example for better understanding. Let $\mathcal{M}$ denotes the $\mathbf{A}_{K}^{mn}$, the space of $m \times n$ matrix and we write the variables $x_{i,j}$, where $1 \le i \le m, 1 \le j \le n$. Let $I^{m,n}_{r}$ be the ideal of $\mathcal{M}$ generated by all $r \times r$ minors of the matrix $(x_{i,j})_{1\le i \le m,1 \le j \le n}$. Now we set $S=K[x_{i,j}]$, then $S/I^{m,n}_{r}$ is the coordinate ring of $\mathscr{L}^{m,n}_{r}$. We fix an integer k, and let $x_{i,j}=x_{i,j}^{(0)}+x_{i,j}^{(1)}t+...+x_{i,j}^{(k)}t^{k}$ for all $(i,j)$. For every polynomial $f(x_{i,j})$ in $I^{m,n}_{r}$, suppose $f(x_{i,j}^{(0)}+x_{i,j}^{(1)}t+...+x_{i,j}^{(k)}t^{k})=f^{(0)}+f^{(1)}t+...+f^{(k)}t^{k}$(mod $t^{k+1}$), then $f^{(0)}, f^{(1)},...,f^{(k)}$ are in the $k$-th jet ideal of $I^{m,n}_{r}$, which we denote by $I^{m,n}_{r,k}$. In fact, one can check that we only need to do the process for the generators of $I^{m,n}_{r}$ to get $I^{m,n}_{r,k}$, i.e. if we assume that $I^{m,n}_{r}=(f_{1},...,f_{s})$, then $I^{m,n}_{r,k}$ is generated by all $f_{i}^{(j)}$, where $1 \le i \le s, 0 \le j \le k$. Now let $R=K[x_{i,j}^{(l)}]$($1 \le i\le m, 1 \le j \le n, 0 \le l \le k$) be the $k$-th jet ring of $S$, then the $k$-th jet ring of $S/I^{m,n}_{r}$ is $R/I^{m,n}_{r,k}$.

\subsection{Basics of Gr\"{o}bner bases}
The theory of Gr\"{o}bner bases is a powerful tool for studying polynomial ideals and forms the basis for many algorithms related to polynomial ideals in computation. In this section, we will fix the ring $S=K[x_{1},...,x_{n}]$ and introduce the basic knowledge of Gr\"{o}bner bases theory. For convenience, we will use the notation of composite indices, that is, given $\alpha=(\alpha_{1},...,\alpha_{n}) \in \mathbb{Z}^{n}_{\ge 0}$, then $x^{\alpha}=x_{1}^{\alpha_{1}}...x_{n}^{\alpha_{n}}$. First, we need to assign a total order to all monomials in $S$, defined specifically as follows:
\begin{definition}
A monomial order $>$ on the set of monomials in $S=K[x_{1},...,x_{n}]$ is a total order satisfying the following conditions:  \\
$(1)$If $\alpha,\beta,\gamma \in \mathbb{Z}^{n}_{\ge 0}$ satisfying $x^{\alpha}>x^{\beta}$, then $x^{\alpha+\gamma}>x^{\beta+\gamma}$. \\
$(2)$For any nonempty subset $A$ of the set of monomials of $S$, $A$ has the minimal element, i.e. there exists $x^{\alpha} \in A$, such that for any $x^{\beta} \in A, x^{\beta} \neq x^{\alpha}$, we have $x^{\beta}>x^{\alpha}$.
\end{definition}

In specific calculations, one monomial order commonly used is graded lexicographic order, where the variables in set $S$ are arranged in descending order as $x_{1}>...>x_{n}$. We provide the graded lexicographic order of $S$ under this arrangement.
\begin{definition}
For any monomials $x^{\alpha},x^{\beta}$ in $S$, we assume $\alpha=(\alpha_{1},...,\alpha_{n}),\beta=(\beta_{1},...,\beta_{n})$, then we define $x^{\alpha}>x^{\beta}$ if and only if $|\alpha|=\alpha_{1}+...+\alpha_{n}>|\beta|=\beta_{1}+...+\beta_{n}$ or $|\alpha|=|\beta|$ and the first non-zero element of $\alpha-\beta$ from left to right is a positive number.
\end{definition}

Similarly we give the reverse graded lexicographic order of $S$:
\begin{definition}
For any monomials $x^{\alpha},x^{\beta}$ in $S$, we assume $\alpha=(\alpha_{1},...,\alpha_{n}),\beta=(\beta_{1},...,\beta_{n})$, then we define $x^{\alpha}>x^{\beta}$ if and only if $|\alpha|=\alpha_{1}+...+\alpha_{n}>|\beta|=\beta_{1}+...+\beta_{n}$ or $|\alpha|=|\beta|$ and the first non-zero element of $\alpha-\beta$ from right to left is a negative number.
\end{definition}

After fixing a monomial order, we can talk about the leading term, leading coefficient and leading monomial of a polynomial, with the following specific definition:
\begin{definition}
Assume $f \in K[x_{1},...,x_{n}]$ is a non-zero polynomial, and suppose $f=\sum_{a_{\alpha} \neq 0} a_{\alpha}x^{\alpha}$, $x^{\alpha_{0}}$ is the maximal element among all the monomials of $f$, then: \\
	$(1)$ The leading term of $f$ is $lt(f):=a_{\alpha_{0}}x^{\alpha_{0}}$. \\
	$(2)$ The leading coefficient of $f$ is $lc(f):=a_{\alpha_{0}}$. \\
	$(3)$ The leading monomial of $f$ is $lm(f):=x^{\alpha_{0}}$.
\end{definition}

Notice that a simple fact is for non-zero polynomials $f,g$, we have $lm(fg)=lm(f)lm(g)$. Therefore, for an ideal $I \subset K[x_{1},...,x_{n}]$, we can also discuss the leading term ideal $LM(I)$, which is defined as the set of all leading monomials of elements in $I$, that is, $LM(I):=\{lm(f)|f \in I\}$. Using $lm(fg)=lm(f)lm(g)$, we know that this set is actually an ideal. It is worth noting that taking only the leading terms of the generating elements of $I$ does not necessarily generate the leading term ideal of $I$. We provide an example from \cite{1} to illustrate this point:
\begin{example}
(\cite{1}, Chapter $2$, Section $5$, Example $2$) Assume $f_{1}=x^{3}-2xy$, $f_{2}=x^{2}y-2y^{2}+x \in K[x,y]$, $I=(f_{1},f_{2})$, and we use the graded lexicographical order in $K[x,y]$. Now because $xf_{2}-yf_{1}=x^{2}$, $x^{2} \in I$, in particular $x^{2} \in LM(I)$, but the leading term of $f_{1}$ is $x^{3}$, and the leading term of $f_{2}$ is $x^{2}y$, so the leading term of $f_{1},f_{2}$ can not generate $LM(I)$.
\end{example}

Next we give the definition and basic properties of Gr\"{o}bner bases.
\begin{definition}
Fix a monomial order on $K[x_{1},...,x_{n}]$, for an ideal $I \subset K[x_{1},...,x_{n}]$, a finite subset $G=\{g_{1},...,g_{s}\} \subset I$ is called a Gr\"{o}bner bases of $I$, if $(lm(g_{1}),...,lm(g_{s}))=LM(I)$. For zero ideal, we define the empty set as its Gr\"{o}bner bases.
\end{definition}

\begin{proposition}
(\cite{1}, Chapter $2$, Section $5$, Corollary $6$) For any ideal $I$ in $K[x_{1},...,x_{n}]$, there exists a Gr\"{o}bner bases $G$, and $G$ is a basis of $I$.
\end{proposition}

In general, computing a given ideal's Gr\"{o}bner bases is not an easy task. In the theory of Gr\"{o}bner bases, its computation mainly relies on the Buchberger criterion. Before explaining this criterion, we need to make some preparations. First, there are propositions and definitions concerning the normal form:
\begin{proposition}
(\cite{1}, Chapter $2$, Section $6$, Proposition $1$) Suppose $I \subset K[x_{1},...,x_{n}]$ is an ideal, $G=\{g_{1},...,g_{s}\} \subset I$, then for any $f \in K[x_{1},...,x_{n}]$, there exists $r \in K[x_{1},...,x_{n}]$ satisfying the following condition: \\
$(1)$ $f=g+r,g \in I$.  \\
$(2)$ Any monomial of $r$ can not be divisible by $lm(g_{1}),...,lm(g_{s})$. \\
We call $r$ a normal form of $f$ with respect to the basis $G$. In particular, when $G$ is a Gr\"{o}bner bases of $I$, $r$ is unique.
\end{proposition}

Please note that in the definition above, the definition of the normal form does not require $G$ to be a Gr\"{o}bner bases. A simple corollary that can be obtained from the proposition above is:
\begin{proposition}
(\cite{1}, Chapter $2$, Section $6$, Corollary $2$) Under the conditions of the proposition above, $f \in I$ if and only if the normal form of $f$ with respect to $G$ is $0$. 
\end{proposition}

Next we give the definition of S-polynomial:
\begin{definition}
Assume $f,g \in K[x_{1},...,x_{n}]$, and $lm(f)=x^{\alpha},lm(g)=x^{\beta}$, $\alpha=(\alpha_{1},...,\alpha_{n}),\beta=(\beta_{1},...,\beta_{n}),\gamma=(\gamma_{1},...,\gamma_{n})$, where $\gamma_{i}=max\{\alpha_{i},\beta_{i}\},i=1,2,...,n$.We set $x^{\gamma}=lcm(x^{\alpha},x^{\beta})$, and we define the S-polynomial of $f$ and $g$ as following:
\begin{flalign*}
S(f,g)=\frac{x^{\gamma}}{lt(f)}f-\frac{x^{\gamma}}{lt(g)}g.
\end{flalign*}
\end{definition}

With the introduction of S-polynomials, we can now present the Buchberger Criterion, which provides a way to determine Gr\"{o}bner bases:
\begin{theorem}
(\cite{1}, Chapter $2$, Section $6$, Theorem $6$) Let $I \subset K[x_{1},...,x_{n}]$ be an ideal, and $G=\{g_{1},...,g_{s}\}$ be a set of basis elements for $I$. Then, $G$ is a Gr\"{o}bner bases for $I$ if and only if, for any 
$1 \le i < j \le s$, $0$ is a normal form of the S-polynomial $S(g_{i},g_{j})$ with respect to $G$.
\end{theorem}

In practical applications, the Buchberger Criterion has many strengthened forms that may be more convenient for our use. One such form replaces the condition of the normal form with the condition that the modulo $G$ reduction is equal to $0$. Below, we provide the specific definition and theorem:
\begin{definition}
Let $G=\{g_{1},...,g_{s}\} \subset K[x_{1},...,x_{n}]$ be a finite set, for any $f \in K[x_{1},...,x_{n}]$, we say $f$ reduces to $0$ modulo $G$, denoted by $f \rightarrow_{G} 0$, if $f$ can be written as $g_{1}h_{1}+...+g_{s}h_{s}$ and for all $i$, when $g_{i}h_{i} \neq 0$, we have $lm(f) \ge lm(g_{i}h_{i})$.
\end{definition}

\begin{theorem}
(\cite{1}, Chapter $2$, Section $9$, Theorem $3$) For an ideal $I \subset K[x_{1},...,x_{n}]$, a basis $G=\{g_{1},...,g_{s}\}$ is a Gr\"{o}bner bases of $I$ if and only if for any $i<j$, $S(g_{i},g_{j}) \rightarrow_{G} 0$.
\end{theorem}

For this new condition $S(g_{i},g_{j}) \rightarrow_{G} 0$, we have convenient methods of determination in some special cases. One common proposition is as follows.
\begin{proposition}
For a finite set $G \subset K[x_{1},...,x_{n}]$, if $f,g \in G$ satisfy that $lm(f)$ and $lm(g)$ are coprime, then $S(f,g) \rightarrow_{G} 0$.
\end{proposition}

\section{the study of $\mathscr{L}^{m,m}_{m,k}$}
In this section, we further study $\mathscr{L}^{m,m}_{m,k}$ and prove Theorem \ref{mtha}. Let $(x_{i,j})_{m \times m}$ be a matrix of size $m \times m$. When computing the jet scheme, we denote the unique defining polynomial of $\mathscr{L}^{m,m}_{m}$ $f(x_{i,j})$. Let $x_{i,j}=x_{i,j}^{(0)}+x_{i,j}^{(1)}t+...+x_{i,j}^{(k)}t^{k}$ and substitute it into $f(x_{i,j})$, denoted as $f(x_{i,j}^{(0)}+x_{i,j}^{(1)}t+...+x_{i,j}^{(k)}t^{k})=f^{(0)}+f^{(1)}t+...+f^{(k)}t^{k}(mod \ t^{k+1})$. Then we have $I^{m,m}_{m,k}=(f^{(0)},...,f^{(k)})$. We now fix the following monomial order: first, let $x_{1,1}^{(k)}>...>x_{1,m}^{(k)}>...>x_{m,1}^{(k)}>...>x_{m,m}^{(k)}>x_{1,1}^{(k-1)}>...>x_{m,m}^{(k-1)}>...>x_{m,m}^{(0)}$. Next, consider the reverse lexicographical order under this ordering. The foundation of our study relies on the following result from \cite{14}, which proves that $f^{(0)},...,f^{(k)}$ form a Gr\"{o}bner bases for $I^{m,m}_{m,k}$ under the above monomial order:

\begin{proposition}
(\cite{14}, Theorem $3.3$) Under the above conditions, $f^{(0)},...,f^{(k)}$ form a Gr\"{o}bner bases for $I^{m,m}_{m,k}$.
\end{proposition}

Having a Gr\"{o}bner bases means that we have the leading term ideal $LM(I^{m,m}_{m,k})$ of $I^{m,m}_{m,k}$. Next, we want to correspond this leading term ideal to an abstract simplicial complex and prove its shellability. To do this, we need to make some preparations. First, we give the definition of an abstract simplicial complex:
\begin{definition}
An abstract simplicial complex $\Delta$ is a family of subsets of a finite set 
$S$, being closed under taking subsets, meaning that for any $X \in \Delta$ and $Y \subset X$, we have $Y \in \Delta$. The elements in $\Delta$ form a partially ordered set under inclusion, and the maximal elements under this partial order are called the facets of $\Delta$. The dimension of a facet $X$ in $\Delta$ is defined as {\color{black}$|X|-1$}, and the dimension of $\Delta$ is defined as the maximal dimension among all facets in $\Delta$.
\end{definition}

\begin{remark}
$(1)$ In studying abstract simplicial complexes, we generally don't focus on the specific elements themselves, but rather on the structure of $\Delta$ as a collection of subsets. Without loss of generality, we can always assume that $S=\{1,2,...,n\}$. \\
$(2)$ When the facets of an abstract simplicial complex are determined, the abstract simplicial complex itself is also determined. Therefore, studying the structure of an abstract simplicial complex is equivalent to characterizing its facets.
\end{remark}

An abstract simplicial complex always corresponds to a ring, which is called the Stanley-Reisner ring, as defined below:
\begin{definition}
Let $S=\{1,2,..,n\}$, and let $\Delta$ be an abstract simplicial complex. Then, the Stanley-Reisner ring of $\Delta$ is defined as $R_{\Delta}=K[x_{1},...,x_{n}]/I_{\Delta}$, where $I_{\Delta}:=(x_{i_{1}}...x_{i_{s}}|\{i_{1},...,i_{s}\} \notin \Delta)$ is called the Stanley-Reisner ideal of $\Delta$.
\end{definition}

For a monomial ideal $I$ in $K[x_{1},...,x_{n}]$, if it is a radical ideal, then the process above can be reversed, meaning that there exists an abstract simplicial complex $\Delta$ such that $I_{\Delta}=I$. In fact, taking the elements of $\Delta$ to be all $\{i_{1},...,i_{s}\}$ such that $x_{i_{1}}...x_{i_{s}} \notin I$ suffices.

Now, returning to our leading term ideal $LM(I^{m,m}_{m,k})$, from the proposition at the beginning of this section, we know that this is a monomial ideal and a radical ideal. Hence, taking the set of all variables $x_{i,j}^{(l)}$ as the total set, there exists an abstract simplicial complex $\Delta$ such that $I_{\Delta}=I^{m,m}_{m,k}$. Next, we provide the definition of shellability and prove that $\Delta$ is shellable:
\begin{definition}
An abstract simplicial complex $\Delta$	is called shellable if all its facets have the same dimension and its facets can be arranged in the following sequence $F_{1},...,F_{e}$, such that for any $1 \le i<j \le e$, there exists an element $v \in F_{j}-F_{i}$, and an index $k<j$, such that $F_{j}-F_{k}=\{v\}$. Such an arrangement is called a shelling of $\Delta$.
\end{definition}

\begin{proposition}
Let $K[x_{i,j}^{(l)},1\le i,j \le m, 0\le l \le k]/LM(I_{m,k}^{m,m})$ be the Stanley-Reisner ring corresponding to $\Delta$, then $\Delta$ is shellable.
\end{proposition}

\begin{proof}
Let $D_0,D_1,...,D_{k}$ denote the sets corresponding to $lm(f^{(0)}),lm(f^{(1)}),...,lm(f^{(k)})$ respectively. Then $\Delta$ contains all sets that do not contain any $D_{i}$, where $i=0,1,...,k$. Note that the intersection of any two sets in $D_0,D_1,...,D_{k}$ is empty. Therefore, if a set $V$ is a facet of $\Delta$, then its complement $V^{c}= \{ v_0,v_1,...,v_{k}\}$, where $v_{i} \in D_{i}$. Now, by fixing a total order on the elements of each $D_{i}$, we can define a total order on the facets of $\Delta$ as follows: for two facets $U$ and $V$, where $U^{c}= \{ u_0,u_1,...,u_{k} \} $ and $V^{c}= \{ v_0,v_1,...,v_{k} \} $, $U<V$ if and only if there exist $0 \le i \le k$ such that $u_{0}=v_{0},...,u_{i-1}=v_{i-1},u_{i}<v_{i}$. Finally, we prove that this ordering yields a shelling of $\Delta$. For any facets $U$ and $V$ satisfying the condition $U<V$, where $u_{0}=v_{0},...,u_{i-1}=v_{i-1},u_{i}<v_{i}$ for some 
$0 \le i \le k$, we have $u_{i} \in V-U$. Let $W$ be a facet such that $W^{c}=\{v_0,...,v_{i-1},u_{i},v_{i+1},...,v_{k}\}$. Then, $W<V$ and $V-W=\{u_{i}\}$, which proves our proposition.
\end{proof}

When an abstract simplicial complex is shellable, the corresponding Stanley-Reisner ring has many desirable properties. Not only can we determine the dimension and Cohen-Macaulay property of this ring, we can also use formulas to calculate its Hilbert series. In this case, it is the Hilbert series of $K[x_{i,j}^{(l)}|1 \le i,j \le m,0 \le l \le k]/LM(I^{m,m}_{m,k})$, which is also the Hilbert series of $K[x_{i,j}^{(l)}|1 \le i,j \le m,0 \le l \le k]/I^{m,m}_{m,k}$. Specifically, we have the following proposition.

\begin{proposition}
(\cite{4}, Theorem $6.3$)Let $\Delta$ be a shellable abstract simplicial complex, $R_{\Delta}$ be its Stanley-Reisner ring, then we have:  \\
$(1)$ $R_{\Delta}$ is Cohen-Macaulay ring. with dimension $1+dim(\Delta)$. \\
$(2)$ Set $d=dim(R_{\Delta})$, $F_{1},...,F_{e}$ is a shelling of $\Delta$, $c(F_{t}):=\{v \in F_{t}| \exists s<t, F_{t}-F_{s}=\{v\} \}$, then the Hilbert series of $R_{\Delta}$ is:
\begin{flalign*}
	\frac{\sum_{j \ge 0}h_{j}z^{j}}{(1-z)^{d}},
\end{flalign*}
where $h_{j}:=|\{t \in \{1,2,...,e\}| \ |c(F_{t})|=j\}|$.
\end{proposition}

By using the proposition above to $K[x_{i,j}^{(l)}|1 \le i,j \le m,0 \le l \le k]/LM(I^{m,m}_{m,k})$, we can immediately get:
\begin{proposition}
$K[x_{i,j}^{(l)}|1 \le i,j \le m,0 \le l \le k]/I^{m,m}_{m,k}$ is a Cohen-Macaulay ring.
\end{proposition}
\begin{proof}
From the proposition above we know $K[x_{i,j}^{(l)}|1 \le i,j \le m,0 \le l \le k]/LM(I^{m,m}_{m,k})$ is Cohen-Macaulay ring, and from Corollary $8.31$ in \cite{15} we know this is equivalent to $K[x_{i,j}^{(l)}|1 \le i,j \le m,0 \le l \le k]/I^{m,m}_{m,k}$ is Cohen-Macaulay ring.
\end{proof}

Proposition 3.7 can also be used to calculate the Hilbert series of $K[x_{i,j}^{(l)}|1 \le i,j \le m,0 \le l \le k]/I^{m,m}_{m,k}$, and thus we obtain Theorem \ref{mtha}.
\begin{theorem}
The Hilbert series of $K[x_{i,j}^{(k)}|1 \le i,j \le m,0 \le l \le k]/I_{m,k}^{m,m}$ is $\frac{(1-z^{m})^{k+1}}{(1-z)^{m^{2}(k+1)}}$($z$ is variable).
\end{theorem}
\begin{proof}
To compute it, we only need to know the value of $|c(V)|$ for each facet $V$. For a given facet $V$, its complement is $V^{c}=\{ v_0,v_1,...v_{k} \}$, and let the number of elements in $D_{i}$ that are less than $V_{i}$ be denoted as $a_{i}$, for $i=0,1,...,k$. Note that an element $v \in c(V)$ if and only if for some $0 \le j \le k$, $v$ is less than $v_{j}$ in $D_{j}$, so $|c(V)|=a_0+a_1+...+a_{k}$. Using the proposition above, we obtain the Hilbert series as follows:
\begin{flalign*}
	\frac{\displaystyle\sum_{j=0}^{\infty}h_{j}z^{j}}{(1-z)^{(m^2-1)(k+1)}} 
	=&\frac{\displaystyle\sum_{j=0}^{\infty}\sum_{\substack{a_0+a_1+...+a_{k}=j \\ 0 \le a_0,...,a_{k} \le m-1}}z^{j}}{(1-z)^{(m^2-1)(k+1)}} 
	=\frac{\displaystyle\sum_{j=0}^{\infty}\sum_{\substack{a_0+a_1+...+a_{k}=j \\ 0 \le a_0,...,a_{k} \le m-1}}z^{a_0+a_1+...+a_{k}}}{(1-z)^{(m^2-1)(k+1)}} \\
	=&\frac{(z^{0}+z^{1}+...z^{m-1})^{k+1}}{(1-z)^{(m^2-1)(k+1)}} 
	=\frac{(1-z^{m})^{k+1}}{(1-z)^{m^{2}(k+1)}}.
\end{flalign*}

\end{proof}

The Hilbert polynomial of $K[x_{i,j}^{(k)}|1 \le i,j \le m,0 \le l \le k]/I_{m,k}^{m,m}$ can also be computed immediately:
\begin{corollary}
The Hilbert polynomial of $K[x_{i,j}^{(k)}|1 \le i,j \le m,0 \le l \le k]/I_{m,k}^{m,m}$ is:
$$P(t)=\displaystyle\sum_{i=0}^{k+1}(-1)^{i}\binom{k+1}{i}\binom{t-mi+m^{2}(k+1)-1}{m^{2}(k+1)-1}.$$ 
\end{corollary}
\begin{proof}
We can compute it directly. Set $N=m(k+1)$,  $(1-z^{m})^{k+1}=\displaystyle\sum_{j=0}^{N}g_{j}z^{j}$, and the specific calculation is as following: 
\begin{flalign*}
	\frac{(1-z^{m})^{k+1}}{(1-z)^{m^{2}(k+1)}}   
	=&\frac{\displaystyle\sum_{j=0}^{N}g_{j}z^{j}}{(1-z)^{m^2(k+1)}}   \\
	=&(\displaystyle\sum_{j=0}^{N}g_{j}z^{j})(\displaystyle\sum_{n=0}^{\infty}\binom{n+m^{2}(k+1)-1}{n}z^{n})   \\
	=&\displaystyle\sum_{j=0}^{N}f_{j}(g_0,...,g_{N})z^{j}+\displaystyle\sum_{j=N+1}^{\infty}(\displaystyle\sum_{l=0}^{N}g_{l}\binom{j-l+m^{2}(k+1)-1}{m^{2}(k+1)-1})z^{j},
\end{flalign*}
where $f_{j}(g_0,...,g_{N})$ denotes the polynomial determined by $g_0,...,g_{N}$. The polynomial we want is exactly $\displaystyle\sum_{l=0}^{N}g_{l}\binom{j-l+m^{2}(k+1)-1}{m^{2}(k+1)-1}$, and we can get the result by calculating each $g_{l}$.
\end{proof}

\section{A Gr\"{o}bner bases of $\mathscr{L}^{3,n}_{3,1}$}
Starting from this section and the subsequent sections, we will study the properties of $\mathscr{L}^{3,n}_{3,1}$. In particular, in this section, we will find a Gr\"{o}bner bases for $\mathscr{L}^{3,n}_{3,1}$. Let $\{x_{i,j}\}(1 \le i \le 3, 1\le j \le n)$ denote the variables of a $3 \times n$ matrix. When computing $\mathscr{L}^{3,n}_{3,1}$, we obtain new variables $\{x_{i,j}^{(0)}\}$ and $\{x_{i,j}^{(1)}\}$. For all $i,j$, we identify $x_{i,j}$ with $x_{i,j}^{(0)}$ and denote $x_{i,j}^{(1)}$ as $y_{i,j}$. First, we assign the following order to all variables:
\begin{flalign*}
y_{1,1}>...>y_{1,n}>y_{2,1}>...>y_{2,n}>y_{3,1}>...>y_{3,n}>x_{1,1}>...>x_{1,n}>x_{2,1}>...>x_{3,n}.
\end{flalign*}
In the subsequent research, we will fix the monomial order as the graded reverse lexicographical order under the above arrangement. Next, in order to define the polynomials for $\mathscr{L}^{3,n}_{3,1}$, we compute the 3x3 minors of the matrix $(x_{i,j}+y_{i,j}t)_{3 \times n}$ and compare the coefficients of the constant term and the term with $t$. This gives us the following two families of polynomials:
\begin{flalign*}
	a_{p,q,r}=det \begin{pmatrix}
		x_{1,p} & x_{1,q} & x_{1,r} \\
		x_{2,p} & x_{2,q} & x_{2,r} \\
		x_{3,p} & x_{3,q} & x_{3,r},
	\end{pmatrix},
\end{flalign*}
\begin{flalign*}
	b_{p,q,r}=det \begin{pmatrix}
		x_{1,p} & x_{1,q} & x_{1,r} \\
		x_{2,p} & x_{2,q} & x_{2,r} \\
		y_{3,p} & y_{3,q} & y_{3,r} 
	\end{pmatrix}+det \begin{pmatrix}
		x_{1,p} & x_{1,q} & x_{1,r} \\
		y_{2,p} & y_{2,q} & y_{2,r} \\
		x_{3,p} & x_{3,q} & x_{3,r} 
	\end{pmatrix}+det \begin{pmatrix}
		y_{1,p} & y_{1,q} & y_{1,r} \\
		x_{2,p} & x_{2,q} & x_{2,r} \\
		x_{3,p} & x_{3,q} & x_{3,r} 
	\end{pmatrix},
\end{flalign*}
where $1 \le p < q < r \le n$. Thus we know all the $a_{p,q,r},b_{p,q,r}$ can define the ideal $I^{3,n}_{3,1}$.

Next we introduce the following three families of polynomials:
\begin{flalign*}
	c_{p,q,r,s}=det \begin{pmatrix}
		x_{1,p} & x_{1,q} & x_{1,r} & x_{1,s} \\
		y_{2,p} & y_{2,q} & y_{2,r} & y_{2,s} \\
		y_{3,p} & y_{3,q} & y_{3,r} & y_{3,s} \\
		x_{3,p} & x_{3,q} & x_{3,r} & x_{3,s}
	\end{pmatrix}+det \begin{pmatrix}
		y_{1,p} & y_{1,q} & y_{1,r} & y_{1,s} \\
		x_{2,p} & x_{2,q} & x_{2,r} & x_{2,s} \\
		y_{3,p} & y_{3,q} & y_{3,r} & y_{3,s} \\
		x_{3,p} & x_{3,q} & x_{3,r} & x_{3,s}
	\end{pmatrix},
\end{flalign*}
\begin{flalign*}
	d_{l,p,q,r,s}=det \begin{pmatrix}
		x_{1,l} & y_{1,p} & y_{1,q} & y_{1,r} & y_{1,s} \\
		x_{2,l} & y_{2,p} & y_{2,q} & y_{2,r} & y_{2,s} \\
		0       & x_{2,p} & x_{2,q} & x_{2,r} & x_{2,s} \\
		x_{3,l} & y_{3,p} & y_{3,q} & y_{3,r} & y_{3,s} \\
		0       & x_{3,p} & x_{3,q} & x_{3,r} & x_{3,s} 
	\end{pmatrix},
\end{flalign*}
\begin{flalign*}
	e_{l,p,q,r,s}=det \begin{pmatrix}
		y_{1,l} & y_{1,p} & y_{1,q} & y_{1,r} & y_{1,s} \\
		y_{2,l} & y_{2,p} & y_{2,q} & y_{2,r} & y_{2,s} \\
		x_{2,l} & x_{2,p} & x_{2,q} & x_{2,r} & x_{2,s} \\
		y_{3,l} & y_{3,p} & y_{3,q} & y_{3,r} & y_{3,s} \\
		x_{3,l} & x_{3,p} & x_{3,q} & x_{3,r} & x_{3,s}
	\end{pmatrix},
\end{flalign*}
where $1 \le p<q<r<s \le n$, and in $d_{l,p,q,r,s}$, $l \le p$, in $e_{l,p,q,r,s}$, $l<p$. We will prove all these $a_{p,q,r},b_{p,q,r},c_{p,q,r,s},d_{l,p,q,r,s},e_{l,p,q,r,s}$ form a Gr\"{o}bner bases of $\mathscr{L}^{3,n}_{3,1}$. Before this, we first prove that these polynomials are indeed in $I^{3,n}_{3,1}$.
\begin{proposition}
$a_{p,q,r},b_{p,q,r},c_{p,q,r,s},d_{l,p,q,r,s},e_{l,p,q,r,s}$ are in $I^{3,n}_{3,1}$.
\end{proposition}
\begin{proof}
We will prove that $c_{p,q,r,s},d_{l,p,q,r,s},e_{l,p,q,r,s}$ are in $I^{3,n}_{3,1}$ respectively. First we calculate $c_{p,q,r,s}$ as following:
\begin{flalign*}
		&c_{p,q,r,s}  \\
		=&det \begin{pmatrix}
			x_{1,p} & x_{1,q} & x_{1,r} & x_{1,s} \\
			y_{2,p} & y_{2,q} & y_{2,r} & y_{2,s} \\
			y_{3,p} & y_{3,q} & y_{3,r} & y_{3,s} \\
			x_{3,p} & x_{3,q} & x_{3,r} & x_{3,s}
		\end{pmatrix}+det \begin{pmatrix}
			y_{1,p} & y_{1,q} & y_{1,r} & y_{1,s} \\
			x_{2,p} & x_{2,q} & x_{2,r} & x_{2,s} \\
			y_{3,p} & y_{3,q} & y_{3,r} & y_{3,s} \\
			x_{3,p} & x_{3,q} & x_{3,r} & x_{3,s}
		\end{pmatrix}+det \begin{pmatrix}
			x_{1,p} & x_{1,q} & x_{1,r} & x_{1,s} \\
			x_{2,p} & x_{2,q} & x_{2,r} & x_{2,s} \\
			y_{3,p} & y_{3,q} & y_{3,r} & y_{3,s} \\
			y_{3,p} & y_{3,q} & y_{3,r} & y_{3,s}
		\end{pmatrix}    \\
		=&y_{3,p}b_{q,r,s}-y_{3,q}b_{p,r,s}+y_{3,r}b_{p,q,s}-y_{3,s}b_{p,q,r}.
\end{flalign*}
This tells us $c_{p,q,r,s}$ is in $I^{3,n}_{3,1}$. Next we compute $e_{l,p,q,r,s}$:
	\begin{flalign*}
		&e_{l,p,q,r,s}   \\
		=&det \begin{pmatrix}
			y_{1,l} & y_{1,p} & y_{1,q} & y_{1,r} & y_{1,s} \\
			y_{2,l} & y_{2,p} & y_{2,q} & y_{2,r} & y_{2,s} \\
			x_{2,l} & x_{2,p} & x_{2,q} & x_{2,r} & x_{2,s} \\
			y_{3,l} & y_{3,p} & y_{3,q} & y_{3,r} & y_{3,s} \\
			x_{3,l} & x_{3,p} & x_{3,q} & x_{3,r} & x_{3,s}
		\end{pmatrix}+det \begin{pmatrix}
			x_{1,l} & x_{1,p} & x_{1,q} & x_{1,r} & x_{1,s} \\
			y_{2,l} & y_{2,p} & y_{2,q} & y_{2,r} & y_{2,s} \\
			y_{2,l} & y_{2,p} & y_{2,q} & y_{2,r} & y_{2,s} \\
			y_{3,l} & y_{3,p} & y_{3,q} & y_{3,r} & y_{3,s} \\
			x_{3,l} & x_{3,p} & x_{3,q} & x_{3,r} & x_{3,s}
		\end{pmatrix}  \\
		=&-y_{2,l}c_{p,q,r,s}+y_{2,p}c_{l,q,r,s}-y_{2,q}c_{l,p,r,s}+y_{2,r}c_{l,p,q,s}-y_{2,s}c_{l,p,q,r}.
	\end{flalign*}
Thus $e_{l,p,q,r,s}$ is also in $I^{3,n}_{3,1}$.
	
Before calculating $d_{l,p,q,r,s}$, we introduce another family of polynomials $c_{p,q,r,s}'$:
$$c_{p,q,r,s}'=det \begin{pmatrix}
		y_{1,p} & y_{1,q} & y_{1,r} & y_{1,s} \\
		y_{2,p} & y_{2,q} & y_{2,r} & y_{2,s} \\
		x_{2,p} & x_{2,q} & x_{2,r} & x_{2,s} \\
		x_{3,p} & x_{3,q} & x_{3,r} & x_{3,s}
	\end{pmatrix}+det \begin{pmatrix}
		x_{1,p} & x_{1,q} & x_{1,r} & x_{1,s} \\
		y_{2,p} & y_{2,q} & y_{2,r} & y_{2,s} \\
		x_{2,p} & x_{2,q} & x_{2,r} & x_{2,s} \\
		y_{3,p} & y_{3,q} & y_{3,r} & y_{3,s}
	\end{pmatrix}.$$
Note that $c_{p,q,r,s}'$ is obtained by replacing the order of the rows and changing the position of the rows, so $c_{p,q,r,s}'$ is also in $I^{3,n}_{3,1}$. Next we calculate $d_{l,p,q,r,s}$:
	\begin{flalign*}
		&d_{l,p,q,r,s}  \\
		=&x_{1,l}det \begin{pmatrix}
			y_{2,p} & y_{2,q} & y_{2,r} & y_{2,s} \\
			x_{2,p} & x_{2,q} & x_{2,r} & x_{2,s} \\
			y_{3,p} & y_{3,q} & y_{3,r} & y_{3,s} \\
			x_{3,p} & x_{3,q} & x_{3,r} & x_{3,s}
		\end{pmatrix}-x_{2,l}det \begin{pmatrix}
			y_{1,p} & y_{1,q} & y_{1,r} & y_{1,s} \\
			x_{2,p} & x_{2,q} & x_{2,r} & x_{2,s} \\
			y_{3,p} & y_{3,q} & y_{3,r} & y_{3,s} \\
			x_{3,p} & x_{3,q} & x_{3,r} & x_{3,s}
		\end{pmatrix} \\
		   &-x_{3,l}det \begin{pmatrix}
			y_{1,p} & y_{1,q} & y_{1,r} & y_{1,s} \\
			y_{2,p} & y_{2,q} & y_{2,r} & y_{2,s} \\
			x_{2,p} & x_{2,q} & x_{2,r} & x_{2,s} \\
			x_{3,p} & x_{3,q} & x_{3,r} & x_{3,s}
		\end{pmatrix}  \\
		=&x_{1,l}det \begin{pmatrix}
			y_{2,p} & y_{2,q} & y_{2,r} & y_{2,s} \\
			x_{2,p} & x_{2,q} & x_{2,r} & x_{2,s} \\
			y_{3,p} & y_{3,q} & y_{3,r} & y_{3,s} \\
			x_{3,p} & x_{3,q} & x_{3,r} & x_{3,s}
		\end{pmatrix}+x_{2,l}det \begin{pmatrix}
			x_{1,p} & x_{1,q} & x_{1,r} & x_{1,s} \\
			y_{2,p} & y_{2,q} & y_{2,r} & y_{2,s} \\
			y_{3,p} & y_{3,q} & y_{3,r} & y_{3,s} \\
			x_{3,p} & x_{3,q} & x_{3,r} & x_{3,s}
		\end{pmatrix} \\
		&+x_{3,l}det \begin{pmatrix}
			x_{1,p} & x_{1,q} & x_{1,r} & x_{1,s} \\
			y_{2,p} & y_{2,q} & y_{2,r} & y_{2,s} \\
			x_{2,p} & x_{2,q} & x_{2,r} & x_{2,s} \\
			y_{3,p} & y_{3,q} & y_{3,r} & y_{3,s}
		\end{pmatrix}-x_{2,l}c_{p,q,r,s}-x_{3,l}c_{p,q,r,s}'  \\
		=&det \begin{pmatrix}
			x_{1,l} & x_{1,p} & x_{1,q} & x_{1,r} & x_{1,s} \\
			0       & y_{2,p} & y_{2,q} & y_{2,r} & y_{2,s} \\
			x_{2,l} & x_{2,p} & x_{2,q} & x_{2,r} & x_{2,s} \\
			0       & y_{3,p} & y_{3,q} & y_{3,r} & y_{3,s} \\
			x_{3,l} & x_{3,p} & x_{3,q} & x_{3,r} & x_{3,s}
		\end{pmatrix}-x_{2,l}c_{p,q,r,s}-x_{3,l}c_{p,q,r,s}'.
	\end{flalign*}
Expanding the determinant of the last step using the second and fourth row, by the Laplace's theorem, it can be written as a linear combination of $a_{i,j,k}$'s, where $i,j,k$ are taken from the set $\{l,p,q,r,s\}$. Therefore, 
$d_{l,p,q,r,s}$ also belongs to $I^{3,n}_{3,1}$. This completes the proof of the proposition.
\end{proof}

Let $\Lambda$ be the set consisting of all $a_{p,q,r},b_{p,q,r},c_{p,q,r,s},d_{l,p,q,r,s},e_{l,p,q,r,s}$. Then, by Theorem $2.16$, to prove that $\Lambda$ forms a Gr\"{o}bner bases of $I^{3,n}_{3,1}$, it suffices to show that for any $\alpha,\beta \in \Lambda$, their S-polynomial $S(\alpha,\beta)$ can be written as $\sum f_{\theta}\theta$, where $\theta \in \Lambda$, and $lm(f_{\theta}\theta) \le lm(S(\alpha,\beta))$. To simplify the problem further, we consider Lemma $2.5$ from \cite{2}:
\begin{lemma}
(\cite{2}, Lemma $2.5$)Suppose $S(\alpha,\beta)=\sum f_{\theta}\theta$, satisfying $\theta \in \Lambda$, $lm(f_{\theta}\theta) \le lm(S(\alpha,\beta))$ and all the variables of $\alpha,\beta$ and $f_{\theta},\theta$ comes from a fixed submatrix $T=\{(i_{s},j_{s}) | i_{1}<...<i_{p},j_{1}<...<j_{q},p,q \in \mathbb{N}\}$, then if replace $T$ with another matrix $T'=\{(i_{s}',j_{s}') | i_{1}'<...<i_{p}',j_{1}'<...<j_{q}',p,q \in \mathbb{N}\}$ such that $T'$ has the same size as $T$, and replace the variables $x_{i_{s},j_{s}},y_{i_{s},j_{s}}$ with $x_{i_{s},j_{s}}',y_{i_{s},j_{s}}'$, we can get another family of polynomials. These polynomials $\alpha',\beta',\theta',f_{\theta}'$ still satisfies $S(\alpha',\beta')=\sum f_{\theta}'\theta'$, $lm(f_{\theta}'\theta') \le lm(S(\alpha',\beta'))$.
\end{lemma}

Note that the variables of the elements in $\Lambda$ all come from a matrix of at most size $3 \times 5$, so for $\alpha,\beta \in \Lambda$, a matrix of size at most $3 \times 10$ is sufficient to contain all the variables. If we can prove that the corresponding $\Lambda$ forms a Gr\"{o}bner bases for the $3 \times 10$ matrix, then by the above lemma, we will have $S(\alpha,\beta)\rightarrow_{\Lambda}0$. According to Proposition $2.17$, when the leading terms of $\alpha$ and $\beta$ are coprime, $S(\alpha,\beta)\rightarrow_{\Lambda}0$ is guaranteed. Therefore, we only need to consider cases where the submatrices of $\alpha$ and $\beta$ overlap, in which case their elements come from a matrix of at most $3 \times 9$. Hence, we only need to consider matrices of size $3 \times 4$ up to $3 \times 9$, and verify that the corresponding $\Lambda$ forms a Gr\"{o}bner bases in each of these cases.

We will use the computational software Singular to verify the cases of matrices of size $3 \times 4$ up to $3 \times 9$. According to the definition of Gr\"{o}bner bases, we only need to calculate $LM(I^{3,n}_{3,1})$ for $n=4,...,9$, and verify that $LM(I^{3,n}_{3,1})$ is exactly generated by the leading terms of the elements in $\Lambda$.

\section{The leading ideal of $I^{3,n}_{3,1}$}
In the previous section, we found a Gr\"{o}bner bases for $I^{3,n}_{3,1}$, so we know that $LM(I^{3,n}_{3,1})$ is generated by the leading terms of the elements in $\Lambda$. We respectively denote the set of leading terms of $a_{p,q,r},b_{p,q,r},c_{p,q,r,s},d_{l,p,q,r,s},e_{l,p,q,r,s}$ as $A,B,C,D,E$, which can be depicted in the following diagram:

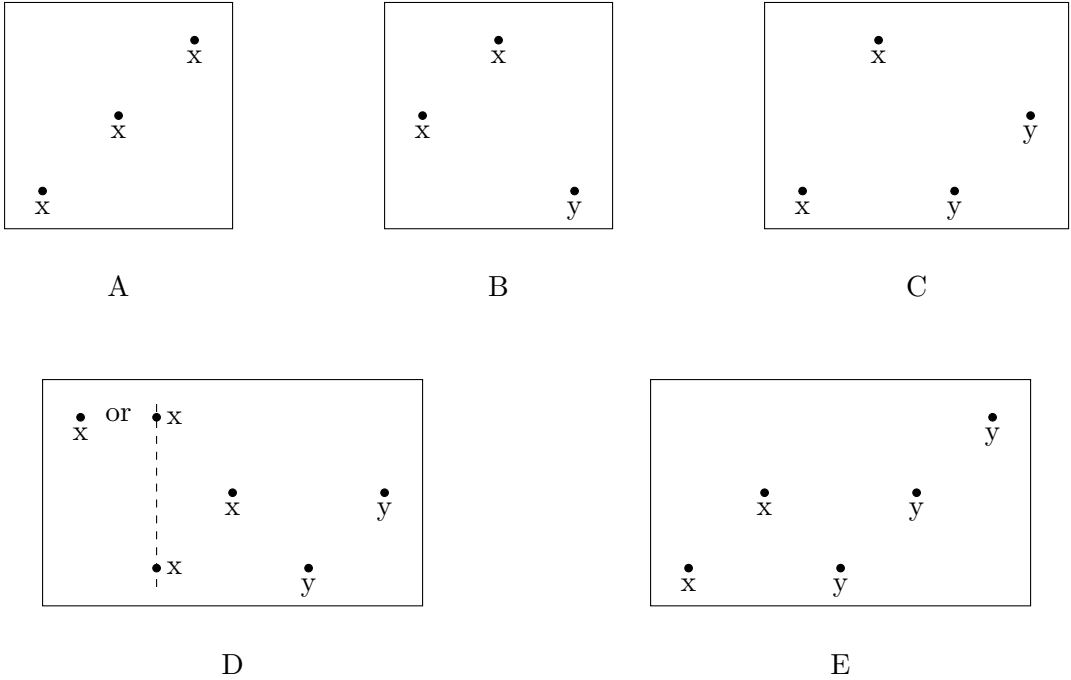
\begin{tikzpicture}
	\draw (0.5,0.5) rectangle (3.5,3.5);
	\draw[fill] (1,1) circle [radius=0.05];
	\node[below] at (1,1) {x};
	\draw[fill] (2,2) circle [radius=0.05];
	\node[below] at (2,2) {x};
	\draw[fill] (3,3) circle [radius=0.05];
	\node[below] at (3,3) {x};
	\node[below] at (2,0) {A};
	
	\draw (5.5,0.5) rectangle (8.5,3.5);
	\draw[fill] (8,1) circle [radius=0.05];
	\node[below] at (8,1) {y};
	\draw[fill] (6,2) circle [radius=0.05];
	\node[below] at (6,2) {x};
	\draw[fill] (7,3) circle [radius=0.05];
	\node[below] at (7,3) {x};
	\node[below] at (7,0) {B};
	
	\draw (10.5,0.5) rectangle (14.5,3.5);
	\draw[fill] (11,1) circle [radius=0.05];
	\node[below] at (11,1) {x};
	\draw[fill] (13,1) circle [radius=0.05];
	\node[below] at (13,1) {y};
	\draw[fill] (14,2) circle [radius=0.05];
	\node[below] at (14,2) {y};
	\draw[fill] (12,3) circle [radius=0.05];
	\node[below] at (12,3) {x};
	\node[below] at (12.5,0) {C};
	
	\draw (1,-4.5) rectangle (6,-1.5);
	\draw[fill] (2.5,-4) circle [radius=0.05];
	\node[right] at (2.5,-4) {x};
	\draw[fill] (4.5,-4) circle [radius=0.05];
	\node[below] at (4.5,-4) {y};
	\draw[fill] (3.5,-3) circle [radius=0.05];
	\node[below] at (3.5,-3) {x};
	\draw[fill] (5.5,-3) circle [radius=0.05];
	\node[below] at (5.5,-3) {y};
	\draw[fill] (1.5,-2) circle [radius=0.05];
	\node[below] at (1.5,-2) {x};
	\draw[fill] (2.5,-2) circle [radius=0.05];
	\node[right] at (2.5,-2) {x};
	\node[below] at(2,-1.75) {or};
	\node[below] at (3.5,-5) {D};
	\draw[dashed] (2.5,-4.25)--(2.5,-1.75);
	
	\draw (9,-4.5) rectangle (14,-1.5);
	\draw[fill] (9.5,-4) circle [radius=0.05];
	\node[below] at (9.5,-4) {x};
	\draw[fill] (11.5,-4) circle [radius=0.05];
	\node[below] at (11.5,-4) {y};
	\draw[fill] (10.5,-3) circle [radius=0.05];
	\node[below] at (10.5,-3) {x};
	\draw[fill] (12.5,-3) circle [radius=0.05];
	\node[below] at (12.5,-3) {y};
	\draw[fill] (13.5,-2) circle [radius=0.05];
	\node[below] at (13.5,-2) {y};
	\node[below] at (11.5,-5) {E};
\end{tikzpicture}

All the sets satisfy the relative position above are contained in these families. 

Since $I^{3,n}_{3,1}$ is a radical ideal, it corresponds to an abstract simplicial complex $\Delta_{LM(I^{3,n}_{3,1})}$, which we will denote as $\Delta_{0}$ for convenience. Next, we aim to characterize the facets of $\Delta_{0}$.

First, based on previous results, we can establish the restrictions of the family $A$ on $\Delta_{0}$. In fact, as Theorem $1$ in \cite{3}, the polynomials defined in the traditional determinant variety itself form a Gr\"{o}bner bases, meaning that when the variables are only $x_{i,j}$, all $a_{p,q,r}$ form a Gr\"{o}bner bases. Therefore, the family $A$ generates the leading term ideal of $I^{3,n}_{3}$. Proposition $6.4$ in \cite{4} also provides a characterization of the corresponding Stanley-Reisner complex it corresponds to:

\begin{proposition}  
(\cite{4}, Proposition $6.4$)For a general $m \times n$ matrix $(x_{i,j})_{m \times n}$, consider the Stanley-Reisner complex $\Delta_{LM(I^{m,n}_{r})}$ corresponding to the ideal $I^{m,n}_{r}$. Then, the facets of $\Delta_{LM(I^{m,n}_{r})}$ correspond to the union of all non-intersecting "paths" respectively from $x_{1,1},x_{2,1},...,x_{r-1}$ to $x_{m,n},...,x_{m,n-r+2}$.
\end{proposition}

Here, the term "path" from $x_{i,j}$ to $x_{p,q}$ refers to a sequence of elements originating from $x_{i,j}$ and ending at $x_{p,q}$, each time moving one step to the right or one step down. The path includes all the elements $x_{u,v}$ visited along the way. An example of a "path" is illustrated in the diagram below.

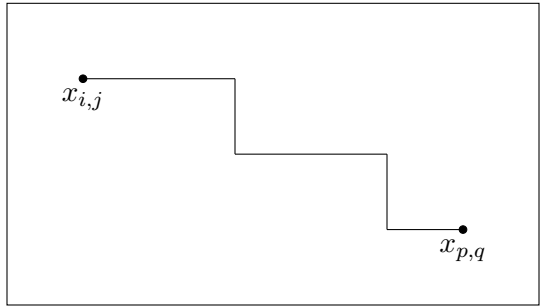
\begin{center}
\begin{tikzpicture}
	\draw (0,0) rectangle (7,4);
	\draw (1,3)--(3,3);
	\draw (3,3)--(3,2);
	\draw (3,2)--(5,2);
	\draw (5,2)--(5,1);
	\draw (5,1)--(6,1);
	\draw[fill] (1,3) circle [radius=0.05];
	\draw[fill] (6,1) circle [radius=0.05];
	\node[below] at (1,3) {$x_{i,j}$};
	\node[below] at (6,1) {$x_{p,q}$};
\end{tikzpicture}
\end{center}

Returning to our original problem, let $F$ be a facet of $\Delta_{0}$, and we consider $F$ and all the elements within $F$ multiplied together to obtain monomials without square factors. Let $F=F_{x}F_{y}$, where $F_{x}$ represents the part composed of $x_{i,j}$ and $F_{y}$ represents the part composed of $y_{i,j}$. Since $F \in \Delta_{0}$, it cannot fully contain any element from the set $A$. In particular, $F_{x}$ cannot completely contain any element from the set $A$. Therefore, when we consider only the variables $x_{i,j}$, $F_{x}$ is contained in some facet of $\Delta_{LM(I^{m,n}_{r})}$. Combined with the previous proposition, we have:

\begin{proposition} 
Let $F=F_{x}F_{y}$ be a facet of $\Delta_{0}$, then $F_{x}$ is a subset of the union of non-intersecting paths from $x_{1,1},x_{2,1}$ to $x_{3,n},x_{3,n-1}$. Specifically, there exists a path $P_{1}$ from $x_{1,1}$ to $x_{3,n}$ and a path $P_{2}$ from $x_{2,1}$ to $x_{3,n-1}$, such that $F_{x} \subset P_{1} \cup P_{2}$ and $P_{1} \cap P_{2}=\phi$. 
\end{proposition}

Next, we consider in detail the specific shape of $P_{1}$ and $P_{2}$. As $P_{1}$
is a path from $x_{1,1}$ to $x_{3,n}$, it must contain exactly two downward movements, and since $P_{1}$ and $P_{2}$ do not intersect, the path $P_{1}$ from the second row to the third row must run from $x_{2,n}$ to $x_{3,n}$. We consider the movement of $P_{1}$ from the first row to the second row, let's denote it as $x_{1,f_{1}}$ to $x_{2,f_{1}}$, then $2 \le f_{1} \le n$. Similarly, $P_{2}$ contains exactly one downward movement, let's denote it from $x_{2,f_{2}}$ to $x_{3,f_{2}}$. Therefore, $P_{1}$ and $P_{2}$ are fully determined by $f_{1}$ and $f_{2}$, and since $P_{1}$ and $P_{2}$ do not intersect, $1 \le f_{2} < f_{1} \le n$. This is illustrated in the following diagram.

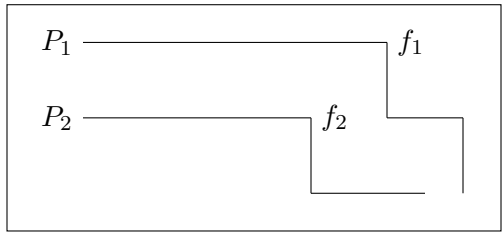
\begin{center}
\begin{tikzpicture}
	\draw (0,0.5) rectangle (6.5,3.5);
	\draw (1,3)--(5,3);
	\draw (5,3)--(5,2);
	\draw (5,2)--(6,2);
	\draw (1,2)--(4,2);
	\draw (4,2)--(4,1);
	\draw (4,1)--(5.5,1);
	\draw (6,2)--(6,1);
	\node[right] at (5,3) {$f_{1}$};
	\node[right] at (4,2) {$f_{2}$};
	\node[left] at (1,3) {$P_{1}$};
	\node[left] at (1,2) {$P_{2}$};
\end{tikzpicture}
\end{center}

We then consider the restrictions of the family $B$ on $F$. If there exist $x_{1,g_{1}},x_{2,g_{2}} \in F_{x}$ such that $g_{2}<g_{1}$, we choose the pair of $(g_{1},g_{2})$ such that $g_{1}$ is the smallest, and among all such pairs we choose $g_{2}$ to be the smallest possible value. In this case, due to the restrictions of family $B$, $y_{3,g_{1}+1},...,y_{3,n} \notin F_{y}$.

We assert that in this situation, $F$ cannot fully contain any elements from families $C$, $D$ or $E$. If $F$ contains the elements $\{x_{3,p},x_{1,q},y_{3,r},y_{2,s}\}(p<q<r<s)$ from family $C$, then $r \le g_{1}$ and $p \ge g_{2}$. In this case, it follows that $q<r \le g_{1},q>p \ge g_{2}$, so $x_{1,q},x_{2,g_{2}} \in F_{x},g_{2}<q<g_{1}$, which contradicts the minimality of $g_{1}$, therefore $F$ does not contain elements in family $C$. If 
$F$ contains the elements $\{x_{1,l},x_{3,p},x_{2,q},y_{3,r},y_{2,s}\}(l \le p<q<r<s)$ in family $D$, then $p \ge g_{2}$, so $x_{2,q}$ can only be in $P_{1}$, thus $q \ge f_{1} \ge g_{1}$. However, $r \le g_{1}$, which is a contradiction, therefore $F$ does not contain elements in family $D$. Similar contradictions can be derived if $F$ contains elements in family $E$, therefore our assertion holds, which tells us that for $F_{x} \subset P_{1} \cup P_{2}$ and $y_{3,g_{1}+1},...,y_{3,n} \notin F_{y}$, $F$ must be in $\Delta_{0}$. Furthermore, since $F$ is a facet of $\Delta_{0}$, we can obtain the following $F_{x},F_{y}$ by adding all elements that satisfy the above conditions as follows:
\begin{flalign*}
	F_{x}=&P_{1} \cup P_{2} - \{x_{1,g_{2}+1},...,x_{1,g_{1}-1},x_{2,1},...,x_{2,g_{2}-1}\},  
\end{flalign*}
\begin{flalign*}
	F_{y}=\{y_{i,j} | 1\le i \le 3, 1 \le j \le n\}-\{y_{3,g_{1}+1},...,y_{3,n}\}.
\end{flalign*}

\begin{tikzpicture}
	\draw (0.5,0.5) rectangle (6.5,3.5);
	\draw (1,3)--(2,3);
	\draw (2,3)--(2,2);
	\draw (2,2)--(4,2);
	\draw (4,2)--(4,1);
	\draw (4,1)--(6,1);
	\draw (3,3)--(5,3);
	\draw (5,3)--(5,2);
	\draw (5,2)--(6,2);
	\draw (6,2)--(6,1);
	\node[below] at (3.5,0) {$x_{i,j}$};
	\node[below] at (2,2) {$g_{2}$};
	\node[left] at (3,3) {$g_{1}$};
	\node[right] at (4,2) {$f_{2}$};
	\node[right] at (5,3) {$f_{1}$};
	\draw (8.5,0.5) rectangle (14.5,3.5);
	\node[below] at (11.5,0) {$y_{i,j}$};
	\draw [dashed] (11,3.25)--(11,0.75);
	\node[left] at (11,3.25) {$g_{1}$};
	\draw (9,3)--(14,3);
	\draw (9,2)--(14,2);
	\draw (9,1)--(11,1);
	\draw[fill] (11,1) circle [radius=0.05];
\end{tikzpicture}

The family of facets shown above is denoted as $\hat{B}$, where the main parameters $g_{1}, g_{2}$ take the value range of $1 \le g_{ 2} <g_{ 1} \le n$. Then we divide the family $F_{ x} $ into two parts, and denote $F_{ B_{ x,2}} =\{ x_{ 1,1},... ,x_{1,g_{2}},x_{2,g_{2}},... ,x_{2,f_{2}},x_{3,f_{2}},... ,x_{3,n}\},F_{B_{x,1}}=\{x_{1,g_{1}},... ,x_{1,f_{1}},x_{2,f_{1}},... ,x_{ 2,n} \}$, where $F_{ B_{ x,2}} $ is the path from $x_{ 1,1} $ to $x_{ 3,n} $, and $F_{ B_{ x,1}} $ is the path from $x_{ 1,g_{ 1}} $ to $x_{ 2,n} $.

When there is no $x_{1,g_{1}},x_{2,g_{2}} \in F_{x}$ satisfying $g_{2}<g_{1}$, then family $B$ no longer restricts $F$. Next we consider family $C$. When there is $x_{1,h_{1}},x_{3,h_{2}} \in F_{x} $ satisfying $h_{1}>h_{2} $, we also take the $h_{1}, h_{2} $ that minimize $h_{ 1} $ and then minimize $h_{2}$ after $h_{1}$ is minimized. Then according to the shape of family $C$, the $y_{i,j}$ in the second row to the third row of the $h_{1}+1$ to the $n$-th column is the path from $y_{2,h_{1}+1}$ to $y_{3,n}$. We assert that such $F$ cannot contain elements of families $D,E$. In fact, if it contains an element of family $D$, denoted as $\{x_{1,l},x_{3,p} x_{2,q},y_{3,r},y_{2,s}\}(l \le p<q<r<s)$, then from $q \ge h_{ 1} $ we know that including $y_{ 3,r},y_{ 2,s} $ will lead to a contradiction. Similarly, including elements of family $E$ also leads to a contradiction. At this point, we can use the maximum property of $F$ to obtain the following $F_{ x},F_{ y}$:
\begin{flalign*}
	F_{x}=P_{1} \cup \{x_{3,h_{2}},...,x_{3,n}\}-\{x_{1,h_{2}+1},...,x_{1,h_{1}-1}\},
\end{flalign*}
\begin{flalign*}
	F_{y}=\{y_{i,j}|1 \le i \le 3, 1 \le j \le h_{1}\} \cup \{y_{1,h_{1}+1},...,y_{1,n}\} \cup F_{C_{y}},
\end{flalign*}

\begin{tikzpicture}
	\draw (0.5,0.5) rectangle (6.5,3.5);
	\draw (1,3)--(3,3);
	\draw (4,3)--(5,3);
	\draw (5,3)--(5,2);
	\draw (5,2)--(6,2);
	\draw (6,2)--(6,1);
	\draw (3,1)--(6,1);
	\draw[fill] (3,1) circle [radius=0.05];
	\draw[fill] (3,3) circle [radius=0.05];
	\draw[dashed] (3,3.25)--(3,0.75);
	\node[left] at (3,1) {$h_{2}$};
	\node[left] at (4,3) {$h_{1}$};
	\node[right] at (5,3) {$f_{1}$};
	\node[below] at (3.5,0) {$x_{i,j}$};
	
	\draw (8.5,0.5) rectangle (14.5,3.5);
	\node[below] at (11.5,0) {$y_{i,j}$};
	\draw [dashed] (12,3.25)--(12,0.75);
	\node[left] at (12,3.25) {$h_{1}$};
	\draw (9,3)--(14,3);
	\draw (9,2)--(13,2);
	\draw (13,2)--(13,1);
	\draw (13,1)--(14,1);
	\draw (9,1)--(12,1);
	\draw[fill] (12,1) circle [radius=0.05];
\end{tikzpicture}

where $F_{ C_{ y}} $ denotes the path from $y_{ 2,h_{ 1} +1} $ to $y_{ 3,n} $. We call the family of facets as shown above $\hat{ C} $, where the main parameters $h_{ 1},h_{ 2} $ take values in the range of $1 \le h_{ 2} <h_{ 1} \le n-1$ (when $h_{ 1} =n$, $F_{y}$ will contain all the $y_{i,j}$, so they will be included in the facet of $h_{ 1} =n-1$). We denote $F_{ C_{ x,2}} =\{ x_{ 3,h_{ 2}},...,x_{ 3,n} \} $ the path from $x_{ 3,h_{ 2}} $ to $x_{ 3,n} $, and $F_{ C_{ x,1}} =\{ x_{ 1,h_{ 1}} ...,x_{ 3,n} \} $ the path from $x_{ 1,h_{ 1}} $ to $x_{ 3,n} $,  respectively. 

If there is no $x_ {1, h_ {1}}, x_ {3, h_ {2}} \in F_ {x} $ satisfy $h_ {1} > h_ {2} $, we finally consider the restriction of family $D,E$. When $F_{x} \cap P_{1}$ has no elements from the first row, $F_{x}$ has no elements from the first row, so we only need to consider the restriction of family $E$ to $F$. Notice $x_ {2, n}, x_ {3, n - 1}$ are not contained in any of $A, B, C, D, E $, so the $x_ {2, n}, x_ {3, n - 1} \in F_ {x} $. Thus there is always $x_{2,k_{1}},x_{3,k_{2}}$ such that $k_{2}<k_{1}$. We take the one of all such pairs that minimizes $k_{1}$, and  minimize $k_{2}$ after guaranteeing $k_{1}$ minimizes. After fixing $k_{1}$, by the shape of family $E$ we know the restriction on $F_{y}$ is that the $y_{i,j}$ from $k_{1}+1$-th column to $n$-th column can only taken to be the non-intersecting paths from $y_{1,k_{1}+1},y_{2,k_{1}+1}$ to $y_{3,n},y_{3,n-1}$. When $k_{1} \le n-2$, combined with the maximality of $F$, we get the following $F_{x},F_{y}$:
\begin{flalign*}
	F_{x}=\{x_{2,1},...,x_{2,k_{2}},x_{3,k_{2}},...,x_{3,n},x_{2,k_{1}},...,x_{2,n}\},
\end{flalign*}
\begin{flalign*}
	F_{y}=\{y_{i,j} | 1 \le i \le 3, 1 \le j \le k_{1}\} \cup F_{E_{y,1}} \cup F_{E_{y,2}},
\end{flalign*}

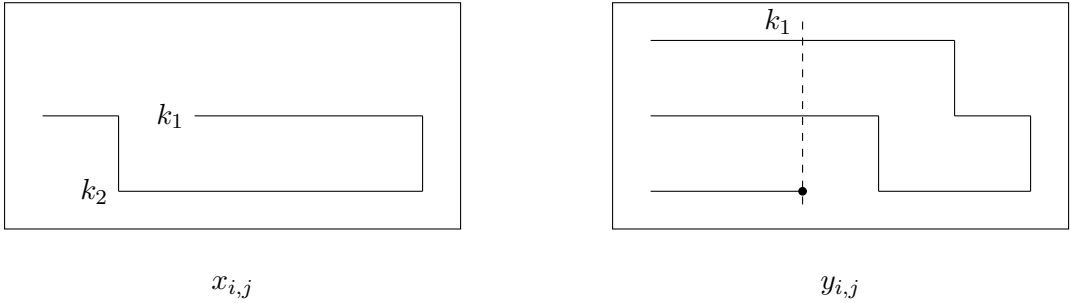
\begin{tikzpicture}
	\draw (0.5,0.5) rectangle (6.5,3.5);
	\draw (1,2)--(2,2);
	\draw (2,2)--(2,1);
	\draw (2,1)--(6,1);
	\draw (3,2)--(6,2);
	\draw (6,2)--(6,1);
	\node[left] at (2,1) {$k_{2}$};
	\node[left] at (3,2) {$k_{1}$};
	\node[below] at (3.5,0) {$x_{i,j}$};
	
	\draw (8.5,0.5) rectangle (14.5,3.5);
	\node[below] at (11.5,0) {$y_{i,j}$};
	\draw [dashed] (11,3.25)--(11,0.75);
	\node[left] at (11,3.25) {$k_{1}$};
	\draw (9,3)--(13,3);
	\draw (13,3)--(13,2);
	\draw (13,2)--(14,2);
	\draw (14,2)--(14,1);
	\draw (9,2)--(12,2);
	\draw (12,2)--(12,1);
	\draw (12,1)--(14,1);
	\draw (9,1)--(11,1);
	\draw[fill] (11,1) circle [radius=0.05];
\end{tikzpicture}

where $F_{E_{y,1}},F_{E_{y,2}}$ denote the disjoint paths of the $k_{1}+1$-th column to the $n$-th column from $y_{1,k_{1}+1},y_{2,k_{1}+1}$ to $y_{3,n},y_{3,n-1}$. We denote the family of facets shown above as $\hat{E}$, where the main parameters $k_{1},k_{2}$ take values in the range $1 \le k_{2}<k_{1} \le n-2$. We write $F_{E_{x,2}}=\{x_{2,1}, ... ,x_{2,k_{2}},x_{3,k_{2}},... ,x_{3,n}\}$ is the path from $x_{2,1}$ to $x_{3,n}$, and $F_{E_{x,1}}=\{x_{2,k_{1}},... ,x_{2,n}\}$ is the path from $x_{2,k_{1}}$ to $x_{2,n}$.

When $k_{1}=n-1,n$, we can also get the same result, but since the $F_{y}$ part is always the total set, it will be contained in the facet when $k_{1}=n-2$ above.

Finally if there is an element of the first row in $F_{x}$, then the same as the case above, combining $x_{2,n},x_{3,n-1} \in F_{x}$, we know that there exists $x_{1,l_{1}},x_{2,l_{2}},x_{3,l_{3}}$ such that $l_{1} \le l_{3}<l_{2}$. We take $l_{2}$ that minimizes $l_{2}$, then $l_{3}$ that minimizes $l_{3}$, then $l_{1}$ that maximizes $l_{1}$. The restriction on $F_{y}$ at this point is that the $y_{i,j}$ from columns $l_{2}+1$ to $n$, and from the second to third rows is the path from $y_{2,l_{2}+1}$ to $y_{3,n}$. We assert that such $F$ can no longer contain any element of the family $E$. Indeed, if $F$ contains $\{x_{3,l},x_{2,p},y_{3,q},y_{2,r},y_{1,s}\}(l<p<q<r<s)$, and since there exists no $x_{1,h_{1}},x_{3,h_{2}} \in F_{x}$ that satisfies $h_{1}>h_{2}$, elements in $F_{x}$ of the first row must satisfy that the second coordinate is less than or equal to $l$, and thus $l_{2} \le p$ by the minimality of $l_{2}$. This contradicts the fact that $F$ contains $y_{3,q},y_{2,r}$. Finally, using the maximality of $F$ yields $F_{x},F_{y}$:

\begin{flalign*}
	F_{x}=\{x_{1,1},...,x_{1,l_{1}},x_{2,l_{1}},...,x_{2,l_{3}},x_{2,l_{2}},...,x_{2,n},x_{3,l_{3}},...,x_{3,n}\},
\end{flalign*}
\begin{flalign*}
	F_{y}=\{y_{i,j}|1 \le i \le 3, 1 \le j \le l_{2}\} \cup \{y_{1,l_{2}+1},...,y_{1,n}\} \cup F_{D_{y}},
\end{flalign*}

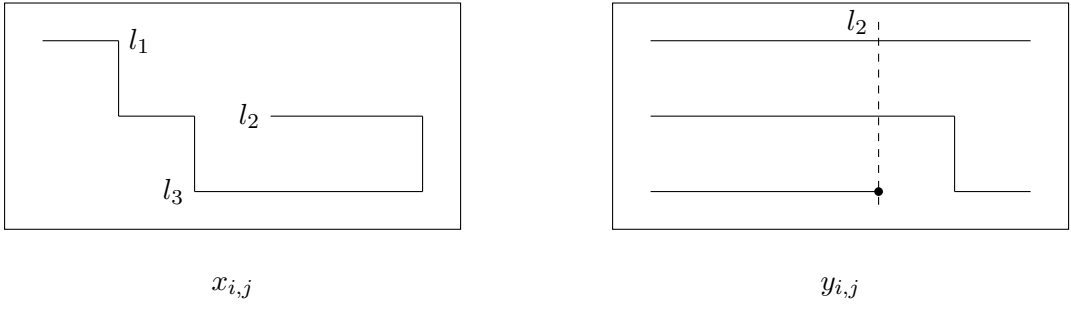
\begin{tikzpicture}
	\draw (0.5,0.5) rectangle (6.5,3.5);
	\draw (1,3)--(2,3);
	\draw (2,3)--(2,2);
	\draw (2,2)--(3,2);
	\draw (3,2)--(3,1);
	\draw (3,1)--(6,1);
	\draw (4,2)--(6,2);
	\draw (6,2)--(6,1);
	\node[left] at (3,1) {$l_{3}$};
	\node[right] at (2,3) {$l_{1}$};
	\node[left] at (4,2) {$l_{2}$};
	\node[below] at (3.5,0) {$x_{i,j}$};
	
	\draw (8.5,0.5) rectangle (14.5,3.5);
	\node[below] at (11.5,0) {$y_{i,j}$};
	\draw [dashed] (12,3.25)--(12,0.75);
	\node[left] at (12,3.25) {$l_{2}$};
	\draw (9,3)--(14,3);
	\draw (9,2)--(13,2);
	\draw (13,2)--(13,1);
	\draw (13,1)--(14,1);
	\draw (9,1)--(12,1);
	\draw[fill] (12,1) circle [radius=0.05];
\end{tikzpicture}

where $F_{D_{y}}$ denotes the path of $y_{i,j}$ from $y_{2,l_{2}+1}$ to $y_{3,n}$. We call the family of facets shown above $\hat{D}$, where the main parameters $l_{2},l_{3}$ take values in the range $1\le l_{3}<l_{2} \le n-1$. We write $F_{D_{x,2}}=\{x_{1,1},...,x_{1,l_{1}},x_{2,l_{1}},...,x_{2,l_{3}},x_{3,l_{3}},... ,x_{3,n}\}$ is the path from $x_{1,1}$ to $x_{3,n}$, $F_{D_{x,1}}=\{x_{2,l_{2}},... ,x_{2,n}\}$ is the path from $x_{2,l_{2}}$ to $x_{2,n}$.

Summarizing, we get all facets of $\Delta_{0}$, four families in total, $\hat{B},\hat{C},\hat{D},\hat{E}$, corresponding to the specific forms shown in the above diagrams as well as the parameters.

\section{A Shelling of $\Delta_{0}$}
In this section we use the facets of $\Delta_{0}$ inscribed in the previous section to prove that $\Delta_{0}$ is shellable. It can be verified that the facets of $\Delta_{0}$ obtained in the previous section all have $4n+4$ elements, and next in order to obtain a shelling of $\Delta_{0}$ we need to give the facets of $\Delta_{0}$ an order.

For the following ordering, we introduce the partial order "$<$" between paths: for two paths $P_{1},P_{2}$ with the same endpoints and starting points in the same row, define $P_{1}<P_{2}$ when the starting point of $P_{1}$ is on the right of $P_{2}$; and when $P_{1},P_ {2}$ have the same starting point, define $P_{1}<P_{2}$ if $P_{1}$ is on the right of $P_{2}$ as shown below. This gives us a partial order relation "$<$", which we extend arbitrarily into a full order relation, still denoted by $<$. For a given path $P$, following the notation in \cite{5}, we call an element $x_{i,j}$ or $y_{i,j}$ "NE turn" if its corresponding $x_{i-1,j},x_{i,j+1}$ or $y_{i-1,j},y_{i,j+2}$ are also in $P$; an element $x_{i,j}$ or $y_{i,j}$ in $P$ is said to be a "SW turn" if the corresponding $x_{i,j-1},x_{i+1,j}$ or $y_{i,j-1},y_{i+1,j}$ are also in $P$.

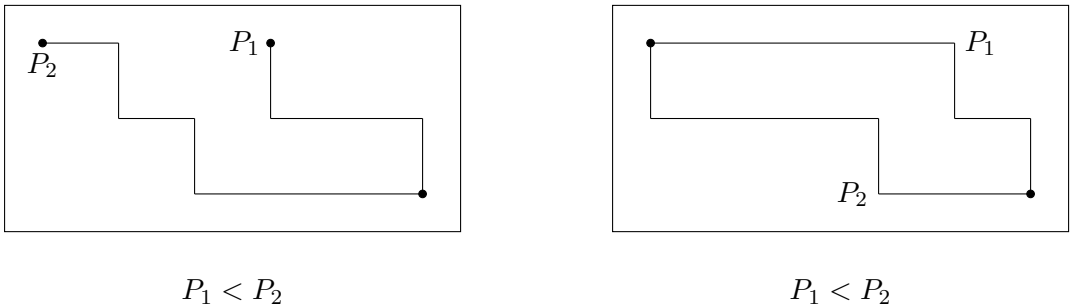
\begin{tikzpicture}
	\draw (0.5,0.5) rectangle (6.5,3.5);
	\draw (1,3)--(2,3);
	\draw (2,3)--(2,2);
	\draw (2,2)--(3,2);
	\draw (3,2)--(3,1);
	\draw (3,1)--(6,1);
	\draw (4,2)--(6,2);
	\draw (6,2)--(6,1);
	\draw (4,3)--(4,2);
	\node[below] at (1,3) {$P_{2}$};
	\node[left] at (4,3) {$P_{1}$};
	\node[below] at (3.5,0) {$P_{1}<P_{2}$};
	\draw[fill] (6,1) circle [radius=0.05];
	\draw[fill] (1,3) circle [radius=0.05];
	\draw[fill] (4,3) circle [radius=0.05];
	
	\draw (8.5,0.5) rectangle (14.5,3.5);
	\node[below] at (11.5,0) {$P_{1}<P_{2}$};
	\node[right] at (13,3) {$P_{1}$};
	\draw (9,3)--(13,3);
	\draw (13,3)--(13,2);
	\draw (13,2)--(14,2);
	\draw (14,2)--(14,1);
	\draw (9,3)--(9,2);
	\draw (9,2)--(12,2);
	\draw (12,2)--(12,1);
	\draw (12,1)--(14,1);
	\node[left] at (12,1) {$P_{2}$};
	\draw[fill] (14,1) circle [radius=0.05];
	\draw[fill] (9,3) circle [radius=0.05];
\end{tikzpicture}

Next we give the ordering rules for the facets of $\Delta_{0}$: for two different facets $P,Q$, if they come from different families, determine their order by $\hat{B}<\hat{D}<\hat{C}<\hat{E}$. If $P,Q$ come from the same family, first compare $P_{x}$ with $Q_{x}$ according to the following rules: first compare the first set of paths, and compare the second set of paths if the first set of paths are the same, e.g., if $P,Q$ both come from $\hat{B}$, then compare $P_{B_{x,1}},Q_{B_{x,1}}$, and when they are the same, then compare $P _{B_{x,2}},Q_{B_{x,2}}$. If $P_{x}$ and $Q_{x}$ are the same, then compare $P_{y}$ and $Q_{y}$: when $P,Q$ belongs to $\hat{C},\hat{D}$, compare $P_{C_{y}},Q_{C_{y}}$ or $P_{D_{y}},Q_{D_{y}}$; when $P,Q$ belongs to $\hat{E}$, we first compare $P_{E_{y,1}},Q_{E_{y,1}}$, and when they are the same, we compare $P_{E_{y,2}},Q_{E_{y,2}}$. In this way we obtain an ordering of the facets of $\Delta_{0}$.

Before we can verify this ordering of the shelling condition, we need the following lemma, which tells us that the "NE turn" is vital for the path:
\begin{lemma}
Let $P_{1},P_{2}$ be two paths with the same start and end points. When $P_{1}$ contains all the "NE turns" of $P_{2}$, $P_{2}$ must be on the right of $P_{1}$.
\end{lemma}
\begin{proof}
Without loss of generality, we assume that $P_{1},P_{2}$ are paths from $x_{1,1}$ to $x_{m,n}$ within the matrix $(x_{i,j})_{m \times n}$. Let $x_{i_{1},j_{1}}, ... ,x_{i_{s},j_{s}}$ be all the "NE turns" of $P_{2}$, then we have $P_{2}=\{x_{1,1}, ... ,x_{1,j_{1}},x_{2,j_{1}},... ,x_{i_{1},j_{1}},... ,x_{i_{s},j_{s}},... ,x_{i_{s},n},... ,x_{m,n}\}$, as shown below:

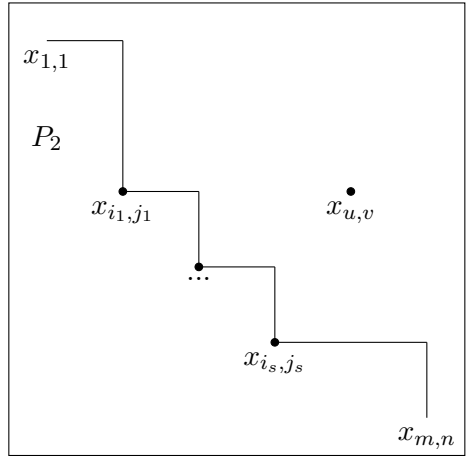
\begin{center}
\begin{tikzpicture}
	\draw (0.5,0.5) rectangle (6.5,6.5);
	\node[below] at (1,5) {$P_{2}$};
	\draw (1,6)--(2,6);
	\node[below] at (1,6) {$x_{1,1}$};
	\draw (2,6)--(2,4);
	\draw (2,4)--(3,4);
	\draw (3,4)--(3,3);
	\draw (3,3)--(4,3);
	\draw (4,2)--(6,2);
	\draw (6,2)--(6,1);
	\node[below] at (6,1) {$x_{m,n}$};
	\draw (4,3)--(4,2);
	\draw[fill] (2,4) circle [radius=0.05];
	\node[below] at (2,4) {$x_{i_{1},j_{1}}$};
	\draw[fill] (3,3) circle [radius=0.05];
	\node[below] at (3,3) {$...$};
	\draw[fill] (4,2) circle [radius=0.05];
	\node[below] at (4,2) {$x_{i_{s},j_{s}}$};
	\draw[fill] (5,4) circle [radius=0.05];
	\node[below] at (5,4) {$x_{u,v}$};
\end{tikzpicture}
\end{center}

Now, using the converse, if $P_{1}$ is not on the left of $P_{2}$, this implies that there exists $x_{u,v} \in P_{1},x_{u,v} \notin P_{2}$ such that there exists $x_{i,j} \in P_{2}$ with $i > u, j < v$, and by the shape of $P_{2}$, there exists $x_{i_{l},j_{ l}}$ such that $i_{l} > u, j_{l} < v$. But $x_{i_{l},j_{l}} \in P_{1}$, which cannot be in the same path as $x_{u,v}$, contradicts. Thus the original proposition holds.
\end{proof}

In the following we verify that this ordering satisfies the conditions of shelling. The plan of verification is: take any facet of $\Delta_{0}$, denoted as $Q$, first find all elements in $c(Q)$, i.e., elements $v \in Q$ such that there exists $R<Q$ with $Q-R=\{v\}$, and then verify that, for any $P<Q$, there exists always some $v$ in $c(Q)$ such that $v \notin P$. For convenience we set each parameter of $Q$ to be the parameter given in the previous section, and then we discuss $Q$ of the following different cases: 

$(1)$ $Q$ is from $\hat{B}$. We assert that $c(Q)$ consists of the following elements: $x_{1,g_{1}} \in c(Q)$ when $x_{1,g_{1}+1} \in Q$; all "NE turns" in $Q_{B_{x,1}},Q_{B_{x,2}}$ except $x_{3,n-1}$.

\begin{tikzpicture}
	\draw (0.5,0.5) rectangle (6.5,3.5);
	\draw (1,3)--(2,3);
	\draw (2,3)--(2,2);
	\draw (2,2)--(4,2);
	\draw (4,2)--(4,1);
	\draw (4,1)--(6,1);
	\draw (3,3)--(5,3);
	\draw (5,3)--(5,2);
	\draw (5,2)--(6,2);
	\draw (6,2)--(6,1);
	\node[below] at (3.5,0) {$x_{i,j}$};
	\node[below] at (2,2) {$g_{2}$};
	\node[left] at (3,3) {$g_{1}$};
	\node[right] at (4,2) {$f_{2}$};
	\node[right] at (5,3) {$f_{1}$};
	\draw (8.5,0.5) rectangle (14.5,3.5);
	\node[below] at (11.5,0) {$y_{i,j}$};
	\draw [dashed] (11,3.25)--(11,0.75);
	\node[left] at (11,3.25) {$g_{1}$};
	\draw (9,3)--(14,3);
	\draw (9,2)--(14,2);
	\draw (9,1)--(11,1);
	\draw[fill] (11,1) circle [radius=0.05];
\end{tikzpicture}

First we verify that the elements above are indeed in $c(Q)$. When $x_{1,g_{1}+1} \in Q$, we can remove $x_{1,g_{1}}$ and replace it with $y_{3,g_{1}+1}$, and verify that this gives us a new facet $R$ and $R<Q$. For the "NE turn" $x_{2,p}$ of $Q_{B_{x,1}}$, we have $x_{1,p},x_{2,p+1} \in Q_{B_{x,1}},x_{1,p+1} \notin Q$, and so replacing $x_{1,p+1}$ with $x_{2,p}$ gives a new facet $R$ and $R<Q$. For a "NE turn" of $Q_{B_{x,2}}$, if $x_{2,q}$ is a "NE turn", then $x_{1,q},x_{2,q+1} \in Q_{B_{x,2}}$. If $x_{1,q+1} \notin Q$, then replacing $x_{2,q}$ with $x_{1,q+1}$ yields a new facet $R$ and $R<Q$; if $x_{1,q+1} \in Q$, then $x_{1,q+1}$ is the starting point of $Q_{B_{x,1}}$ and $x_{1,q+2} \in Q_{B_ {x,1}}$, so removing $x_{1,q}$ and adding $y_{1,q+2}$ gives the new facet $R$ and $R<Q$. If $x_{3,r}$ is a "NE turn" and $r<n-1$, then $x_{2,r},x_{3,r+1} \in Q_{B_{x,2}}$. When $x_{2,r+1} \notin Q$, it is sufficient to replace $x_{3,r}$ with $x_{2,r+1}$; when $x_{2,r+1} \in Q_{B_{x,1}}$, $x_{1,r+1},x_{2,r+2} \in Q_{B_{x,1}}$, so $x_{1,r+2} \notin Q$. Replacing $x_{3,r}$ with $x_{1,r+2}$ gives the new facet $R$ and $R<Q$.

We then show that none of the remaining elements in $Q$ are in $c(Q)$. For $v \in c(Q)$, note that all facets in $\hat{B}$ contain $x_{1,1},x_{3,n},x_{2,n}$, so $v$ cannot be taken as them. Since $R<Q$, $v \in Q_{x}$. When $v$ is not an element in the assertion, $v$ cannot be taken to be the starting or ending point of $Q_{B_{x,1}},Q_{B_{x,2}}$. Let $v_{1},v_{2}$ be the points before and after $v$ on the same path as $v$, respectively. Below we first consider the difference between the two paths in $R_{x}$ and $Q_{x}$.

When $R_{B_{x,1}}$ starts at the same point as $Q_{B_{x,1}}$, $Q_{y}=R_{y}$ and the number of elements of $Q_{x}$ is the same as that of $R_{x}$. If $v \in Q_{B_{x,1}}$, we assert that $Q_{B_{x,1}}-\{v\} \subset R_{B_{x,1}}$. Indeed, for $u \in Q_{B_{x,1}}-\{v\}$, when $u$ is in the first row, consider the point pair $(u,x_{2,g_{2}})$ we know that $u \in R_{B_{x,1}}$; when $u$ is in the second row, consider the point pair $(u,x_{3,f_{2}})$ we know that $u \in R_{B_{x,1}}$. Since $v$ is not a "NE turn", $v$ can only be a "SW turn" or the midpoint of a long segment. If $v$ is the middle point, then by $v_{1},v_{2} \in R_{B_{x,1}}$ we know that $v \in R_{B_{x,1}}$, a contradiction; if $v$ is the "SW turn", then $R_{B_{x,1}}$ is the "SW turn" of $v$. The changing of a point of "SW turn" to "NE turn" makes $R_{B_{x,1}}>Q_{B_{x,1}}$, a contradiction. If $v \in Q_{B_{x,2}}$, similarly we assert that $Q_{B_{x,1}}-\{v\} \subset R_{B_{x,1}}$. Indeed, for $u \in Q_{B_{x,2}}-\{v\}$, when $u$ is in the third row, consider the point pair $(u,x_{2,n})$ we know that $u \in R_{B_{x,2}}$; and when $u$ is in the second row, consider the point pair $(u,x_{1,f_{1}})$ we know that $u \in R_{B_{x,2}}$ and when $u$ is in the first row, consider $(u,x_{1,g_{1}})$ we know $u \in R_{B_{x,2}}$ since $x_{1,g_{1}}$ is still the starting point of $R_{B_{x,1}}$. Therefore $v_{1},v_{2} \in R_{B_{x,2}}$, so again if $v$ is the midpoint, then $v \in R_{B_{x,2}}$, a contradiction; if $v$ is the "SW turn", then $R_{B_{x,2}}$ is the point where one of the "SW turns" of $Q_{B_{x,2}}$ is placed. A "SW turn" is replaced by a "NE turn", and $Q_{B_{x,1}}=R_{B_{x,1}}$, so that $R>Q$, a contradiction.

When $R_{B_{x,1}}$ starts at a different point than $Q_{B_{x,1}}$, since $R<Q$, $R-Q=\{y_{1,g_{1}+1}\}$ and $R_{B_{x,1}}$ starts at $x_{1,g_{1}+1}$. At this point, since $v$ is not $x_{1,g_{1}}$, $x_{1,g_{1}} \in R_{B_{x,2}}$ and thus $x_{1,g_{1}-1},x_{2,g_{1}} \in R_{B_{x,2}}$. This shows that $x_{1,g_{1}-1},x_{2,g_{1}} \in Q_{B_{x,2}}$, so $x_{2,g_{1}-1} \in Q_{B_{x,2}}$, but $x_{2,g_{1}-1} \in R$, so $v=x_{2,g_{1}-1}$, when $v$ is "NE turn", a contradiction. In summary we state all the elements in $c(Q)$.

We then prove that all facets $P$ satisfying $P<Q$ cannot contain $c(Q)$.

Indeed, if there exists a facet $P$ that contains $c(Q)$ and satisfies $P<Q$, we first prove that the start point of $P_{B_{x,1}}$ is the same as $Q_{B_{x,1}}$. If $x_{1,g_{1}} \notin c(Q)$, when $g_{1}=n$, since $P<Q$, we have $P_{B_{x,1}}=Q_{B_{x,1}}$; and when $g_{1}<n$, $x_{2,g_{1}}$ is the "NE turn", thus $x_{2,g_{1}} \in P$. We assert that $x_{2,g_{1}} \in P_{B_{x,1}}$. In fact, notice that there exists a "NE turn" of $Q_{B_{x,2}}$ in the third row, whose column coordinates are less than $g_{1}$, and if it is $x_{3,n-1}$, then $g_{1}=n$, a contradiction, and hence it is in $c(Q)$, and since it is in a different path than $x_{2,g_{1}}$, it is in $c(Q)$. In summary $x_{2,g_{1}} \in P_{B_{x,1}}$. Combined with $P<Q$, we know that $P_{B_{x,1}}$ starts at the same point as $Q_{B_{x,1}}$. If $x_{1,g_{1}} \in c(Q)$, notice that there is a  "NE turn" in $Q$ with second coordinates less than $g_{1}$, and since $P$ contains $c(Q)$, this element and $x_{1,g_{1}}$ can not be in the same path, but only $x_{1,g_{1}} \in P_{B_{x,1}}$. Since $P<Q$, the starting point of $P_{B_{x,1}}$ is also $x_{1,g_{1}}$. This tells us that $P_{B_{x,1}}$ and $Q_{B_{x,1}}$ are two paths with the same starting and ending points. For a "NE turn" $v$ in $Q_{B_{x,1}}$, there must exist a "NE turn" in the third row of $Q_{B_{x,2}}$ with the second coordinate less than the second coordinate of $v$, and hence $v \in P_{B_{x,1}}$. By Lemma $6.1$ we know that $Q_{B_{x,1}}$ is on the right of or coincides with $P_{B_{x,1}}$, but $P<Q$, so $Q_{B_{x,1}}=P_{B_{x,1}}$. Also since $P_{B_{x,2}}$ contains all the "NE turns" of $Q_{B_{x,2}}$, $Q_{B_{x,2}}$ is on the right of $P_{B_{x,2}}$ by Lemma $6.1$, a contradiction. Thus a facet $P$ satisfying $P<Q$ cannot contain $c(Q)$.

$(2)$ $Q$ is from $\hat{D}$. We assert that $c(Q)$ consists of the following elements: the "NE turn" of $Q_{D_{y}}$ and $Q_{D_{x,2}}$; when $Q_{D_{y}}$ has no "NE turn", $y_{3,n} \in c(Q)$; $x_{2,l_{2}} \in c(Q)$.

\begin{tikzpicture}
	\draw (0.5,0.5) rectangle (6.5,3.5);
	\draw (1,3)--(2,3);
	\draw (2,3)--(2,2);
	\draw (2,2)--(3,2);
	\draw (3,2)--(3,1);
	\draw (3,1)--(6,1);
	\draw (4,2)--(6,2);
	\draw (6,2)--(6,1);
	\node[left] at (3,1) {$l_{3}$};
	\node[right] at (2,3) {$l_{1}$};
	\node[left] at (4,2) {$l_{2}$};
	\node[below] at (3.5,0) {$x_{i,j}$};
	
	\draw (8.5,0.5) rectangle (14.5,3.5);
	\node[below] at (11.5,0) {$y_{i,j}$};
	\draw [dashed] (12,3.25)--(12,0.75);
	\node[left] at (12,3.25) {$l_{2}$};
	\draw (9,3)--(14,3);
	\draw (9,2)--(13,2);
	\draw (13,2)--(13,1);
	\draw (13,1)--(14,1);
	\draw (9,1)--(12,1);
	\draw[fill] (12,1) circle [radius=0.05];
\end{tikzpicture}

We start by stating that the elements above are indeed in $c(Q)$. First is $x_{2,l_{2}}$, when $l_{2}=n-1$, replacing $x_{2,n-1}$ with $x_{1,n}$ gives the new facet $R$ and $R<Q$; when $l_{2}<n-1$, exactly one of $y_{3,l_{2}+1},y_{2,l_{2}+2}$ is not in $Q$, replacing it with $x_{2,l_{2}}$ is sufficient. For the "NE turn" $y_{3,p}$ of $Q_{D_{y}}$, it is sufficient to replace it with $y_{2,p+1}$. When $Q_{D_{y}}$ has no "NE turn", it means that $Q_{D_{y}}=\{y_{2,l_{2}+1}, ... ,y_{2,n},y_{3,n}\}$, in which case it is sufficient to replace $y_{3,n}$ with $x_{1,l_{2}}$. For the "NE turn" of $Q_{D_{x,2}}$, if it is $x_{2,q}$, it is sufficient to replace it with $x_{1,q+1}$; if it is $x_{3,r}$, it is sufficient to replace it with $x_{2,r+1}$ when $x_{2,r+1} \notin Q$, and when $x_{2,r+1}$ is notin $Q$, it is sufficient to replace it with $x_{1,q+1}$. $_{2,r+1} \in Q$, it means that $x_{2,r+1}=x_{2,l_{2}}$, in this case, if $l_{2}<n-1$, replace $x_{3,r}$ with the one that $y_{3,l_{2}+1},y_{2,l_{2}+2}$ is not in $Q$, and if $l_{2}=n-1$, replace $x_{1,n}$ with $x_{3,r}$ satisfies the condition.

We then show that the remaining elements are not in $c(Q)$. Now suppose that there is $v$ that is not an element of the above and is in $c(Q)$. Let $Q-R=\{v\}$, and we discuss the different cases of $R$. When $R \in \hat{B}$, if $Q_{D_{y}}$ has a "NE turn", then there exists $y_{2,j}$ which is not in $Q$, and hence the only element of $R$ which is not in $Q$ is $y_{2,j}$. But $Q_{x}$ cannot be $R_{x}$ by decreasing it by one element or by not moving it, a contradiction, so there is no "NE turn" of $Q_{D_{y}}$. This tells us that $Q_{D_{y}}=\{y_{2,l_{2}+1},... ,y_{2,n},y_{3,n}\}$, and since $v$ is not $y_{3,n}$, it follows that $y_{3,n} \in R_{y}$, which means that $R_{B_{x,1}}$ starts at $x_{1,n}$. Thus there can only be $l_{2}=n-1$ and $v=x_{2,l_{2}}$, a contradiction.

When $R \in \hat{D}$, if $v \in Q_{y}$, then the number of elements of $R_{x}$ is the same as or one more than $Q_{x}$, but $R<Q$, so the number of elements of $R_{x}$ is the same as that of $Q_{x}$, and in fact $R_{x}=Q_{x}$. This means that $R_{D_{y}}$ and $Q_{D_{y}}$ have the same starting and ending points and differ by only one element, which is the same as the case in $(1)$ where $v$ can only be a "NE turn", a contradiction. If $v \in Q_{x}$, since $v$ is not the starting or ending point of $Q_{D_{x,1}},Q_{D_{x,2}}$, let $v_{1},v_{2}$ be the two neighboring points of $v$. When $v \in Q_{D_{x,1}}$, for $u \in Q_{D_{x,1}}-\{v\}$, consider the point pair $(u,x_{3,l_{3}})$ we know that $u \in R_{D_{x,1}}$. Therefore $v_{1},v_{2} \in R_{D_{x,1}}$, and $v$ can only be the midpoint of the line segment, so $v \in R_{D_{x,1}}$, a contradiction. When $v \in Q_{D_{x,2}}$, if the starting point of $R_{D_{x,1}}$ is the same as that of $Q_{D_{x,1}}$, we assert that $Q_{D_{x,2}}-\{v\} \subset R_{D_{x,2}}$. Indeed, for $u \in Q_{D_{x,2}}-\{v\}$, when $u$ is in the first row, $u$ can only be in $R_{D_{x,2}}$; when $u$ is in the third row, considering the point pair $(u,x_{2,l_{2}})$ it is known that $u \in Q_{D_{x,2}}$; and when $u$ is in the second row, since $x_{2, l_{2}}$ is the starting point of $R_{D_{x,1}}$, so $u$ can only be in $R_{D_{x,2}}$. This tells us that $v_{1},v_{2} \in R_{D_{x,2}}$, so the same proof as in $(1)$ shows that $v$ can only be a "NE turn". If $R_{D_{x,1}}$ starts at a different point than $Q_{D_{x,1}}$, then $x_{2,l_{2}} \in R_{D_{x,2}}$ and there are more elements of $R$ than $Q$ in $R_{y}$, considering the total number of elements, we know that $x_{2,l_{2}+1}$ is the starting point of $R_{D_{x,1}}$. This tells us that $x_{2,l_{2}-1},x_{3,l_{2}} \in R_{D_{x,2}}$, so $v=x_{3,l_{2}-1}$ is the "NE turn", a contradiction. The rest of the elements of $v=x_{3,l_{2}-1}$ are not in $c(Q)$.

Then we prove that all facets $P$ satisfying $P<Q$ cannot contain $c(Q)$.

If there is a $P$ that satisfies the above condition, when $P \in \hat{B}$, first consider $c(Q) \cap Q_{D_{y}}$: the coordinates of the columns of the elements in $c(Q) \cap Q_{D_{y}}$ are greater than $l_{2}$, and therefore the coordinates of the columns of the starting point of $P_{B_{x,1}}$ are greater than $l_{2}$. This tells us that $x_{2,l_{2}} \in P_{B_{x,2}}$, but $x_{3,l_{3}}$ also belongs to $c(Q) \subset P$ as the "NE turn" of $Q_{D_{x,2}}$, and $x_{3,l_{3}}$ can only be in $P_{B_{x,2}}$, a contradiction. When $P \in \hat{D}$, $x_{2,l_{2}},x_{3,l_{3}} \in P$ necessarily belongs to two different paths, and hence there is $x_{2,l_{2}} \in P_{D_{{x,1}}}$. And since $P<Q$, $x_{2,l_{2}}$ is the starting point of $P_{D_{x,1}}$. Thus $P_{D_{x,2}}$ contains all the "NE turns" of $Q_{D_{x,2}}$, and by Lemma $6.1$ $Q$ is on the right of $P$, contradicting $P<Q$.

$(3)$ $Q$ is from $\hat{C}$. We assert that $c(Q)$ consists of the following elements: $x_{1,h_{1}},x_{3,h_{2}} \in c(Q)$; the "NE turn" of $Q_{C_{x,1}},Q_{C_{y}}$; when $Q_{C_{y}}$ does not have "NE turn", $y_{3,n}\in c(Q)$.

\begin{tikzpicture}
	\draw (0.5,0.5) rectangle (6.5,3.5);
	\draw (1,3)--(3,3);
	\draw (4,3)--(5,3);
	\draw (5,3)--(5,2);
	\draw (5,2)--(6,2);
	\draw (6,2)--(6,1);
	\draw (3,1)--(6,1);
	\draw[fill] (3,1) circle [radius=0.05];
	\draw[fill] (3,3) circle [radius=0.05];
	\draw[dashed] (3,3.25)--(3,0.75);
	\node[left] at (3,1) {$h_{2}$};
	\node[left] at (4,3) {$h_{1}$};
	\node[right] at (5,3) {$f_{1}$};
	\node[below] at (3.5,0) {$x_{i,j}$};
	
	\draw (8.5,0.5) rectangle (14.5,3.5);
	\node[below] at (11.5,0) {$y_{i,j}$};
	\draw [dashed] (12,3.25)--(12,0.75);
	\node[left] at (12,3.25) {$h_{1}$};
	\draw (9,3)--(14,3);
	\draw (9,2)--(13,2);
	\draw (13,2)--(13,1);
	\draw (13,1)--(14,1);
	\draw (9,1)--(12,1);
	\draw[fill] (12,1) circle [radius=0.05];
\end{tikzpicture}

First we verify that these elements are indeed in $c(Q)$. As in the previous case, the "NE turn" of $Q_{C_{x,1}},Q_{C_{y}}$ is naturally in $c(Q)$. When $Q_{C_{y}}$ has no "NE turn", $Q_{C_{y}}=\{y_{2,h_{1}+1},... ,y_{2,n},y_{3,n}\}$, then replacing $y_{3,n}$ with $x_{2,h_{2}}$ yields a new facet $R$ with $R<Q$. For $x_{1,h_{1}}$, when $x_{1,h_{1}+1} \in Q$, if $h_{1}<n-1$, then exactly one of $y_{3,h_{1}+1},y_{2,h_{1}+2}$ is not in $Q_{y}$, and replacing $x_{1,h_{1}}$ with it satisfies the condition; if $h_{1}=n-1 $, then $Q_{y}$ is all $y_{i,j}$, when replacing $x_{1,h_{1}}$ with $x_{2,h_{2}}$ is sufficient. When $x_{1,h_{1}+1} \notin Q$, it is sufficient to replace $x_{1,h_{1}}$ with $x_{2,h_{2}}$. For $x_{3,h_{2}}$, if $x_{1,h_{2}+1} \notin Q$, it is sufficient to replace $x_{3,h_{2}}$ with it. If $x_{1,h_{2}+1} \in Q$, we have $h_{2}+1=h_{1}$. When $x_{1,h_{1}+1} \in Q$, if $h_{1}<n-1$, then exactly one of $y_{3,h_{1}+1},y_{2,h_{1}+2}$ is not in $Q_{y}$, and the condition can be satisfied by replacing $x_{3,h_{2}}$ with it; and if $h_{1}=n-1$, then $Q_{y}$ is all $y_ {i,j}$, at which point it is sufficient to replace $x_{3,h_{2}}$ with $x_{2,h_{2}}$. When $x_{1,h_{1}+1} \notin Q$, if $h_{1}<n-1$, then exactly one of $y_{3,h_{1}+1},y_{2,h_{1}+2}$ is not in $Q_{y}$, and the condition can be satisfied by replacing $x_{3,h_{2}}$ with it; and if $h_{1}=n-1$, then $Q_{y}$ for all $y_{i,j}$, at which point it is sufficient to replace $x_{3,h_{2}}$ with $x_{1,n}$.

We then show that the remaining elements are not in $c(Q)$. Now suppose that there is $v$ that is not an element of the above and is in $c(Q)$, $Q-R=\{v\}$, and we discuss the different cases of $R$. When $R \in \hat{B}$, if $Q_{D_{y}}$ has a "NE turn", then there exists $y_{2,j}$ which is not in $Q$, and hence the only element of $R$ which is not in $Q$ can only be $y_{2,j}$. But $Q_{x}$ cannot be $R_{x}$ by decreasing it by one element or by not moving it, a contradiction, so there is no "NE turn" of $Q_{D_{y}}$. This tells us that $Q_{D_{y}}=\{y_{2,l_{2}+1},... ,y_{2,n},y_{3,n}\}$, and since $v$ is not $y_{3,n}$, it follows that $y_{3,n} \in R_{y}$, which means that $R_{B_{x,1}}$ starts at $x_{1,n}$. Thus there can only be $h_{1}=n-1$, $v=x_{1,h_{1}}$ or $v=x_{3,h_{2}}$, a contradiction.

When $R \in \hat{D}$, $x_{3,h_{2}} \in R_{x}$ since $v$ is not $x_{3,h_{2}}$, notice that there is no element of the first row of $R_{D_{x,1}}$, and therefore none of the elements of the first row of $Q_{C_{x,1}}$ are in $R$, so $v=x_{1,h_{1}}$, a contradiction.

When $R \in \hat{C}$, if $v \in Q_{y}$, then the number of elements of $R_{x}$ is the same as or one more than $Q_{x}$, but $R<Q$, so the number of elements of $R_{x}$ is the same as that of $Q_{x}$, and in fact $R_{x}=Q_{x}$. This means that $R_{C_{y}}$ and $Q_{C_{y}}$ have the same starting and ending points and differ by only one element, which is the same as the case in $(1)$ where $v$ can only be a "NE turn", a contradiction. If $v \in Q_{x}$, since $v$ is not $x_{3,h_{2}},x_{1,h_{1}}$, they are still in $R_{x}$, which tells us that the starting point of $R_{C_{x,2}}$ has column coordinates less than or equal to $h_{2}$, and so $x_{1,h_{1}}$ is still the starting point of $R_{C_{x,1}}$. When the $R_{C_{x,2}}$ starting point is not $x_{3,h_{2}}$, $x_{1,h_{2}} \notin R$, so $v=x_{1,h_{2}}$ and can only be replaced up to $x_{3,h_{2}-1}$, which is a contradiction, so the $R_{C_{x,2}}$ starting point is $x_{3,h_{2}}$. At this point, only $R_{C_{x,1}}$ can be changed, and $R_{C_{x,1}}$ has the same starting and ending points, so $v$ can only be the "NE turn", a contradiction. In summary we show that the elements in $c(Q)$.

Then we prove that all facets $P$ satisfying $P<Q$ cannot contain $c(Q)$.

If there is a $P$ satisfying the above condition when $P \in \hat{B}$, first consider $c(Q) \cap Q_{B_{y}}$: the coordinates of the columns of the elements in $c(Q) \cap Q_{B_{y}}$ are greater than $h_{1}$, and therefore the coordinates of the columns of the starting point of $P_{B_{x,1}}$ are greater than $h_{1}$. And because $x_{1,h_{1}},x_{3,h_{2}} \in c(Q) \subset P$, this tells us that $x_{1,h_{1}} \in P_{B_{x,1}}$, a contradiction. When $P \in \hat{D}$, since $x_{1,h_{1}},x_{3,h_{2}}$ can only be in two different paths, $x_{1,h_{1}} \in P_{D_{{x,2}}}$, a contradiction. When $P \in \hat{C}$, the coordinates of the starting column of $P_{C_{x,2}}$ are less than or equal to $h_{2}$ since $x_{3,h_{2}} \in P$, and thus $x_{1,h_{1}} \in P_{C_{x,1}}$, which tells us that the coordinates of the starting column of $P_{C_{x,1}}$ is less than or equal to $h_{1} $. And since $P<Q$, the $P_{C_{x,1}},P_{C_{x,2}}$ starting point is the same as $Q_{C_{x,1}},Q_{C_{x,2}}$. And since the "NE turn" of $Q_{C_{x,1}},Q_{C_{y}}$ is in $P$, it follows from Lemma 6.1 that $Q$ is on the right of $P$, which is a contradiction of $P<Q$.

$(4)$ $Q$ is from $\hat{E}$. We assert that $c(Q)$ consists of the following elements: the "NE turn" of $Q_{E_{x,2}},Q_{E_{y,1}},Q_{E_{y,2}}$; when $Q_{E_{y,1}}$ has no "NE turn", $y_{2,n} \in c(Q)$; $y_{3,n-1} \in c(Q)$ when $Q_{E_{y,2}}$ has no "NE turn"; $x_{2,k_{1}} \in c(Q)$.

\begin{tikzpicture}
	\draw (0.5,0.5) rectangle (6.5,3.5);
	\draw (1,2)--(2,2);
	\draw (2,2)--(2,1);
	\draw (2,1)--(6,1);
	\draw (3,2)--(6,2);
	\draw (6,2)--(6,1);
	\node[left] at (2,1) {$k_{2}$};
	\node[left] at (3,2) {$k_{1}$};
	\node[below] at (3.5,0) {$x_{i,j}$};
	
	\draw (8.5,0.5) rectangle (14.5,3.5);
	\node[below] at (11.5,0) {$y_{i,j}$};
	\draw [dashed] (11,3.25)--(11,0.75);
	\node[left] at (11,3.25) {$k_{1}$};
	\draw (9,3)--(13,3);
	\draw (13,3)--(13,2);
	\draw (13,2)--(14,2);
	\draw (14,2)--(14,1);
	\draw (9,2)--(12,2);
	\draw (12,2)--(12,1);
	\draw (12,1)--(14,1);
	\draw (9,1)--(11,1);
	\draw[fill] (11,1) circle [radius=0.05];
\end{tikzpicture}

First we verify that these elements are indeed in $c(Q)$. As in the previous case, we know that the "NE turn" of $Q_{E_{y,1}},Q_{E_{y,2}}$ are in $c(Q)$. For the "NE turn" $x_{3,k_{2}}$ of $Q_{E_{x,2}}$, if $x_{2,k_{2}+1} \notin Q$, it is sufficient to replace $x_{3,k_{2}}$ with it; if $x_{2,k_{2}+1} \in Q$, it is $Q _{E_{x,1}}$ starting point, and if $k_{1}<n-2$, then exactly one of $y_{3,k_{1}+1},y_{2,k_{1}+2},y_{1,k_{1}+3}$ is not in $Q$, and it is sufficient to replace $x_{2,k_{2}}$ by it; and if $k_{1}=n-2$, then $Q_{y}$ is all $y_ {i,j}$, at which point it is sufficient to replace $x_{2,k_{2}}$ with $x_{1,1}$. When $Q_{E_{y,1}}$ has no "NE turn", $Q_{E_{y,1}}=\{y_{1,k_{1}+1},... ,y_{1,n},y_{2,n},y_{3,n}\}$, when it is sufficient to replace $y_{2,n}$ with $x_{1,1}$. When $Q_{E_{y,2}}$ does not have a "NE turn", $Q_{E_{y,2}}=\{y_{2,k_{1}+1},... ,y_{2,n-1},y_{3,n-1}\}$, and it is sufficient to replace $y_{3,n-1}$ with $x_{1,1}$.

We then show that the remaining elements are not in $c(Q)$. Now suppose that there is $v$ which is not an element of the above and is in $c(Q)$, and we discuss the different cases of $R$. When $R$ is not in $\hat{E}$, $x_{1,1} \in R$, so $R-Q=\{x_{1,1}\}$. Since $R_{y}$ contains all the elements of the first row, $Q_{y}$ contains $y_{1,n}$, thus $Q_{E_{y,1}}=\{y_{1,k_{1}+1},... ,y_{1,n},y_{2,n},y_{3,n}\}$. Since the elements $x_{i,j}$ in $\hat{B},\hat{C}$ all have at least two elements in the first row, $R$ cannot be in them, and thus $R$ can only be in $\hat{D}$. Since $v$ is not $y_{2,n}$, $y_{2,n} \in R$, so $R_{y}$ contains all elements of the second row in $y_{i,j}$. And because $v$ is not $y_{3,n-1}$, $y_{3,n-1} \in R$, which tells us that $R_{D_{x,1}}$ starts at $x_{2,n-1}$. But $v$ is not $x_{3,k_{2}},x_{2,k_{1}}$, so they belong to $R$, but they are not in the same path, which shows that the coordinates of the columns of the starting point of $R_{D_{x,1}}$ are at most $k_{1}$, a contradiction.

When $R \in \hat{E}$, if $v \in Q_{y}$, then the number of elements of $R_{x}$ is the same as or one more than $Q_{x}$, but $R<Q$, so the number of elements of $R_{x}$ is the same as that of $Q_{x}$, and in fact $R_{x}=Q_{x}$. This means that $R_{E_{y,1}},R_{E_{y,2}}$ and $Q_{E_{y,1}},Q_{E_{y,2}}$ have the same starting and ending points and differ by only one element, the same as in $(1)$, which shows that $v$ can only be a "NE turn", a contradiction, so $v \in Q_{x}$. Since $v$ is not $x_{3,k_{2}},x_{2,k_{1}}$, they belong to $R$, but they are not in the same path, which means that the coordinates of the starting column of $R_{E_{x,1}}$ are at most $k_{1}$, and since $R<Q$, $R_{E_{x,1}}=Q_{E_{x,1}}$, $Q_{E_{x,2}}$ only changes by one element compared to $R_{E_{x,2}}$, a contradiction to the previous proof that $v$ can only be a "NE turn". In summary we have the elements of $c(Q)$.

Next we prove that all facets $P$ satisfying $P<Q$ cannot contain $c(Q)$.

If there is a $P$ satisfying the above condition, when $P \in \hat{B}$, first consider the part of $y_{i,j}$: notice that the "NE turn" of $Q_{E_{y,2}}$, or $x_{3,n-1}$, are the elements in third row with column coordinates greater than $k_{1}$, so the starting column coordinates of $P _{B_{x,1}}$ has starting column coordinates greater than $k_{1}$. Also, since $x_{2,k_{1}} \in c(Q) \subset P$, $x_{2,k_{1}} \in P_{B_{x,2}}$, but $x_{3,k_{2}}$ can only be in $P_{B_{x,2}}$ as well, contradiction. When $P \in \hat{D}$, notice that both elements in $c(Q) \cap Q_{y}$ are in the second and third rows and are not on the same path, so that the coordinates of the starting column of $P_{D_{x,1}}$ are greater than $k_{1}$, and thus $x_{2,k_{1}} \notin P_{D_{x,1}}$. This tells us that $x_{2,k_{1}},x_{3,k_{2}}$ are all in $P_{D_{x,2}}$, a contradiction. When $P \in \hat{C}$, the same can be shown that the column coordinates of the starting point of $P_{C_{x,1}}$ is greater than $k_{1}$, which contradicts $x_{2,k_{1}} \in P$. When $P \in \hat{E}$, since $x_{2,k_{1}},x_{3,k_{2}}$ can't be in the same path, $x_{2,k_{1}} \in P_{E_{{x,1}}},x_{3,k_{2}} \in P_{E_{{x,2}}}$, which tells us that $Q_{x}$ is on the right side of $P_{x}$,  combining with $P<Q$ we know that $P_{x}=Q_{x}$. And finally by Lemma $6.1$ we know that $Q_{y}$ is to the right of $P_{y}$, a contradiction.

\section{The calculation of Hilbert series of $I^{3,n}_{3,1}$}
In this section, we will use the formulas given in the previous section and the specific form of $c(Q)$ for each facet $Q$ of $\Delta_{0}$ to compute the Hilbert series of $I^{3,n}_{3,1}$. Before that, we need a technical lemma from \cite{5}:
\begin{lemma}
(\cite{5}, Proposition $3$) $(1)$ There are integers $a_{1} \le a_{2} \le ... \le a_{s}, a_{1}' \ge ... \ge a_{s}'$, $b_{1} \le b_{2} \le ... \le b_{s}, b_{1}' \ge ... \ge b_{s}'$, then the number of disjoint paths from the lattice points $(a_{1},a_{1}'),... ,(a_{s},a_{s}')$, respectively, to $(b_{1},b_{1}'),... ,(b_{s},b_{s}')$ is:
\begin{flalign*}
	det(\binom{b_{j}-a_{i}+b_{j}'-a_{i}'}{b_{j}-a_{i}}_{1 \le i,j \le s}).
\end{flalign*}
$(2)$ Suppose we have integers $a_{1} \le a_{2} \le ... \le a_{s}, a_{1}' > ... > a_{s}'$, $b_{1} < b_{2} < ... < b_{s}, b_{1}' \ge ... \ge b_{s}'$, then the number of disjoint paths from the lattice points $(a_{1},a_{1}'),...,(a_{s},a_{s}')$, respectively, to $(b_{1},b_{1}'),...,(b_{s},b_{s}')$ with exactly $r$ "NE turns" is:
\begin{flalign*}
	\sum_{r_{1}+... +r_{s}=r}det(\binom{b_{j}-a_{i}+i-j}{r_{i}+i-j}\binom{b_{j}'-a_{i}'-i+j}{r_{i}}_{1 \le i,j \le s}).
\end{flalign*}
\end{lemma}

In the following we will, for each facet $Q$, compute the number of elements of $c(Q)$.

$(1)$ $Q \in \hat{B}$. First recall the elements of $c(Q)$ given in the previous section: $x_{1,g_{1}} \in c(Q)$ when $x_{1,g_{1}+1} \in Q$; all "NE turns" in $Q_{B_{x,1}},Q_{B_{x,2}}$ except $x_{3,n-1}$. For a determined $Q$, we can additionally add "$0$-th row" to the $x_{i,j}$ part and take $Q_{B_{x,1}}'=\{x_{0,2}, ... ,x_{0,g_{1}}\} \cup Q_{B_{x,1}}$ is the path from $x_{0,2}$ to $x_{2,n}$, and $Q_{B_{x,2}}'=Q_{B_{x,2}}-\{x_{3,n}\}$ is the path from $x_{1,1}$ to $x_{3,n-1}$, then the element in $c(Q)$ at that time is the path from $x_{1,1}$ to $x_{3,n-1}$, then the elements of $c(Q)$ are all the "NE turns" of $Q_{B_{x,1}}',Q_{B_{x,2}}'$, as shown in the figure.

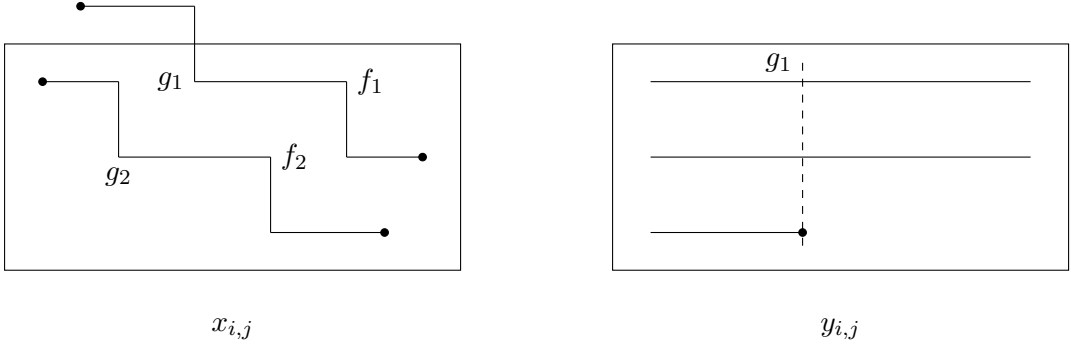
\begin{tikzpicture}
	\draw (0.5,0.5) rectangle (6.5,3.5);
	\draw (1.5,4)--(3,4);
	\draw (3,4)--(3,3);
	\draw (1,3)--(2,3);
	\draw (2,3)--(2,2);
	\draw (2,2)--(4,2);
	\draw (4,2)--(4,1);
	\draw (4,1)--(5.5,1);
	\draw (3,3)--(5,3);
	\draw (5,3)--(5,2);
	\draw (5,2)--(6,2);
	\node[below] at (3.5,0) {$x_{i,j}$};
	\node[below] at (2,2) {$g_{2}$};
	\node[left] at (3,3) {$g_{1}$};
	\node[right] at (4,2) {$f_{2}$};
	\node[right] at (5,3) {$f_{1}$};
	\draw (8.5,0.5) rectangle (14.5,3.5);
	\node[below] at (11.5,0) {$y_{i,j}$};
	\draw [dashed] (11,3.25)--(11,0.75);
	\node[left] at (11,3.25) {$g_{1}$};
	\draw (9,3)--(14,3);
	\draw (9,2)--(14,2);
	\draw (9,1)--(11,1);
	\draw[fill] (11,1) circle [radius=0.05];
	\draw[fill] (1,3) circle [radius=0.05];
	\draw[fill] (1.5,4) circle [radius=0.05];
	\draw[fill] (6,2) circle [radius=0.05];
	\draw[fill] (5.5,1) circle [radius=0.05];
\end{tikzpicture}

By the lemma given at the beginning of this section, the number of $Q$ satisfying $c(Q)=r$ is:
\begin{flalign*}
	\sum_{r_{1}+r_{2}=r}\binom{n-2}{r_{1}}\binom{n-2}{r_{2}}(\binom{2}{r_{1}}\binom{2}{r_{2}}-\binom{2}{r_{1}-1}\binom{2}{r_{2}+1}).
\end{flalign*}
Therefore, there are one $Q$ satisfying $c(Q)=0$, $2n-4$ many $Q$ satisfying $c(Q)=1$, $(n-2)(2n-3)$ many $Q$ satisfying $c(Q)=2$, $frac{2(n-1)(n-2)(n-3)}{3}$ many $Q$ satisfying $c(Q)=3$, $\frac{(n-1)(n-2)^{2}(n-3)}{12}$ many $Q$ satisfying $c(Q)=4$.

$(2)$ $Q \in \hat{D}$. First recall the composition of the elements of $c(Q)$ given in the previous section: the "NE turn" of $Q_{D_{x,2}},Q_{D_{y}}$; when there is no "NE turn" of $Q_{D_{y}}$, $y_{3,n } \in c(Q)$; $x_{2,l_{2}} \in c(Q)$. Noting that there is exactly one element of $c(Q)$ in $Q_{y}$, and the fact that $x_{2,l_{2}} \in c(Q)$ always holds, it is only necessary to consider the number of "NE turns" of $Q_{D_{x,2}}$. When the parameters $l_{1},l_{2},l_{3}$ are determined, $Q_{x}$ is determined, and there are a total of $n-l_{2}$ possibilities for $Q_{y}$, so it is sufficient to consider only the parameters $l_{1},l_{2},l_{3}$, as follows.

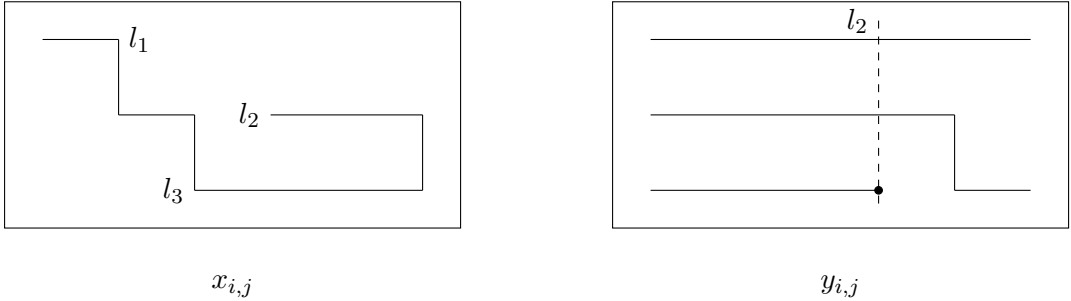
\begin{tikzpicture}
	\draw (0.5,0.5) rectangle (6.5,3.5);
	\draw (1,3)--(2,3);
	\draw (2,3)--(2,2);
	\draw (2,2)--(3,2);
	\draw (3,2)--(3,1);
	\draw (3,1)--(6,1);
	\draw (4,2)--(6,2);
	\draw (6,2)--(6,1);
	\node[left] at (3,1) {$l_{3}$};
	\node[right] at (2,3) {$l_{1}$};
	\node[left] at (4,2) {$l_{2}$};
	\node[below] at (3.5,0) {$x_{i,j}$};
	
	\draw (8.5,0.5) rectangle (14.5,3.5);
	\node[below] at (11.5,0) {$y_{i,j}$};
	\draw [dashed] (12,3.25)--(12,0.75);
	\node[left] at (12,3.25) {$l_{2}$};
	\draw (9,3)--(14,3);
	\draw (9,2)--(13,2);
	\draw (13,2)--(13,1);
	\draw (13,1)--(14,1);
	\draw (9,1)--(12,1);
	\draw[fill] (12,1) circle [radius=0.05];
\end{tikzpicture}

$(2.1)$ $|c(Q)|=3$. At this point $Q_{D_{x,2}}$ has only one "NE turn", so $l_{1}=l_{3}<l_{2} \le n-1$, and thus there are $\sum_{l_{2}=2}^{n-1}(l_{2}-1)(n -l_{2})=\frac{n(n-1)(n-2)}{6}$ many $Q$.

$(2.2)$ $|c(Q)|=4$. At this point $Q_{D_{x,2}}$ has two "NE turns", so $l_{1}<l_{3}<l_{2} \le n-1$, and thus there are a total of $\sum_{l_{2}=3}^{n-1}\frac{(l_{2} -1)(l_{2}-2)(n-l_{2})}{2}=\frac{n(n-1)(n-2)(n-3)}{24}$ many $Q$.

$(3)$ $Q \in \hat{C}$. First recall that $c(Q)$ consists of the following elements: $x_{1,h_{1}},x_{3,h_{2}} \in c(Q)$; the "NE turn" of $Q_{C_{x,1}},Q_{C_{y}}$; when $Q_{C_{y}}$ does not have "NE turn", $y_{3,n}\in c(Q)$. Since $Q_{y}$ has exactly one element in $c(Q)$ and $x_{1,h_{1}},x_{3,h_{2}} \in c(Q)$ always holds, the number of elements in $c(Q)$ depends on the number of "NE turns" in $Q_{C_{x,1}}$. Notice that when the parameters $h_{1},h_{2},f_{1}$ are determined, $Q_{x}$ is determined, and there are $n-h_{1}$ choices for $Q_{y}$, so it is sufficient to consider only the above three parameters, as follows.

\begin{tikzpicture}
	\draw (0.5,0.5) rectangle (6.5,3.5);
	\draw (1,3)--(3,3);
	\draw (4,3)--(5,3);
	\draw (5,3)--(5,2);
	\draw (5,2)--(6,2);
	\draw (6,2)--(6,1);
	\draw (3,1)--(6,1);
	\draw[fill] (3,1) circle [radius=0.05];
	\draw[fill] (3,3) circle [radius=0.05];
	\draw[dashed] (3,3.25)--(3,0.75);
	\node[left] at (3,1) {$h_{2}$};
	\node[left] at (4,3) {$h_{1}$};
	\node[right] at (5,3) {$f_{1}$};
	\node[below] at (3.5,0) {$x_{i,j}$};
	
	\draw (8.5,0.5) rectangle (14.5,3.5);
	\node[below] at (11.5,0) {$y_{i,j}$};
	\draw [dashed] (12,3.25)--(12,0.75);
	\node[left] at (12,3.25) {$h_{1}$};
	\draw (9,3)--(14,3);
	\draw (9,2)--(13,2);
	\draw (13,2)--(13,1);
	\draw (13,1)--(14,1);
	\draw (9,1)--(12,1);
	\draw[fill] (12,1) circle [radius=0.05];
\end{tikzpicture}

$(3.1)$ $|c(Q)|=3$. At this point $Q_{C_{x,1}}$ has no "NE turn", so $f_{1}=n,h_{2}<h_{1} \le n-1$, so there are $\sum_{h_{1}=2}^{n-1}(h_{1}-1)(n-h_{1})=\frac{n(n-1)(n-2)}{6}$ many $Q$. 

$(3.2)$ $|c(Q)|=4$. At this point $Q_{C_{x,1}}$ has a "NE turn", so $h_{2}<h_{1} \le f_{1} \le n-1$, so there are a total of $\sum_{h_{1}=2}^{n-1}(h_{1}-1)(n-h_{1})^{2} =\frac{n(n-1)^{2}(n-2)}{12}$ many $Q$.

$(4)$ $Q \in \hat{E}$. First recall that $c(Q)$ consists of the following elements: the "NE turn" of $Q_{E_{x,2}},Q_{E_{y,1}},Q_{E_{y,2}}$; when $Q_{E_{y,1}}$ has no "NE turn", $y_{2,n} \in c(Q)$; $y_{3,n-1} \in c(Q)$ when $Q_{E_{y,2}}$ has no "NE turn"; $x_{2,k_{1}} \in c(Q)$. Similarly $c(Q) \cap Q_{y}$ always has two elements, $x_{2,k_{1}} \in c(Q)$ always holds, and $Q_{E_{x,2}}$ has exactly one "NE turn", so that $c(Q)$ always has $4$ elements. Since $Q_{x}$ is determined when $k_{1},k_{2}$ is determined, by the lemma at the beginning of this section, there are $\frac{(n-k_{1})(n-k_{1}-1)}{2}$ ways to take $Q_{y}$, and so we only need to consider the parameters $k_{1},k_{2}$ as shown in the figure.

\begin{tikzpicture}
	\draw (0.5,0.5) rectangle (6.5,3.5);
	\draw (1,2)--(2,2);
	\draw (2,2)--(2,1);
	\draw (2,1)--(6,1);
	\draw (3,2)--(6,2);
	\draw (6,2)--(6,1);
	\node[left] at (2,1) {$k_{2}$};
	\node[left] at (3,2) {$k_{1}$};
	\node[below] at (3.5,0) {$x_{i,j}$};
	
	\draw (8.5,0.5) rectangle (14.5,3.5);
	\node[below] at (11.5,0) {$y_{i,j}$};
	\draw [dashed] (11,3.25)--(11,0.75);
	\node[left] at (11,3.25) {$k_{1}$};
	\draw (9,3)--(13,3);
	\draw (13,3)--(13,2);
	\draw (13,2)--(14,2);
	\draw (14,2)--(14,1);
	\draw (9,2)--(12,2);
	\draw (12,2)--(12,1);
	\draw (12,1)--(14,1);
	\draw (9,1)--(11,1);
	\draw[fill] (11,1) circle [radius=0.05];
\end{tikzpicture}

Since $k_{1},k_{2}$ satisfies $k_{2}<k_{1} \le n-2$, the number of facets is $\sum_{k_{1}=2}^{n-2}\frac{(n-k_{1})(n-k_{1}-1)(k_{1}-1)}{2}=\frac{n(n-1)(n-2)(n-3)}{24}$.

Summarizing, we get:
\begin{flalign*}
	h_{0}=&1,h_{1}=2n-4,h_{2}=(n-2)(2n-3), \\
	h_{3}=&\frac{2(n-1)(n-2)(n-3)}{3}+\frac{n(n-1)(n-2)}{6}+\frac{n(n-1)(n-2)}{6} \\
	=&(n-1)(n-2)^{2}, \\
	h_{4}=&\frac{(n-1)(n-2)^{2}(n-3)}{12}+\frac{n(n-1)(n-2)(n-3)}{24}+\frac{n(n-1)^{2}(n-2)}{12}+\frac{n(n-1)(n-2)(n-3)}{24} \\
	=&\frac{(n-1)^{2}(n-2)^{2}}{4}.
\end{flalign*}
The Hilbert series of $K[x_{i,j},y_{i,j}|1 \le i \le 3, 1 \le j \le n]/I^{3,n}_{3,1}$ (with $z$ being the variable) is obtained from Proposition $3.7$, which is Theorem \ref{mthb}.
\begin{flalign*}
	&\frac{1+(2n-4)z+(n-2)(2n-3)z^{2}+(n-1)(n-2)^{2}z^{3}+\frac{(n-1)^{2}(n-2)^{2}}{4}z^{4}}{(1-z)^{4n+4}}   \\
	=&(\frac{1+(n-2)z+\frac{(n-1)(n-2)}{2}z^{2}}{(1-z)^{2n+2}})^{2}.
\end{flalign*}

For the traditional determinantal varieties, i.e., the bottom space $K[x_{i,j}|1 \le i \le 3, 1 \le j \le n]/I^{3,n}_{3}$, the Hilbert series can be given by Theorem $6.9$ in \cite{4} as exactly $\frac{1+(n-2)z+\frac{(n-1)(n-2)}{2 }z^{2}}{(1-z)^{2n+2}}$, i.e., the Hilbert series of $\mathscr{L}^{3,n}_{3,1}$ is the square of the Hilbert series of $\mathscr{L}^{3,n}_{3}$.

\section{The irreducible components of $\mathscr{L}^{m,n}_{r,k}$}
The irreducible components of $\mathscr{L}_{r,k}^{m,n}$ is an interesting topic to study. In \cite{12}, the authors have given the number of irreducible components of a general $\mathscr{L}_{r,k}^{m,n}$ and have characterized them using the $GL_{m} \times GL_{n}$ action. In this section, we will first introduce the theory, and then give the specific polynomials to define these components. To introduce the theory, we need to do some preparations.

Let $G$ denotes the group $GL_{m} \times GL_{n}$, and we have the following natural action of $G$ on $\mathcal{M}$, the space of $m \times n$ matrix:
\begin{flalign*}
	G \times &\mathcal{M} \longrightarrow \mathcal{M}, \\
	(g,h) \times &A \mapsto gAh^{-1}.
\end{flalign*}
Applying the theory of jet schemes and arc spaces, we can get the action of the corresponding jet schemes and arc spaces:
\begin{flalign*}
	G_{n} \times &\mathcal{M}_{n} \longrightarrow \mathcal{M}_{n}, \\
	G_{\infty} \times &\mathcal{M}_{\infty} \longrightarrow \mathcal{M}_{\infty}.
\end{flalign*}
The authors of \cite{12} have proved that, every orbit of this action has a standard form, with the following proposition:
\begin{proposition}
	(Proposition $3.2$ in \cite{12}) 
	Let $\lambda =(\lambda_{1},...,\lambda_{m})$ with $\infty \ge \lambda_{1} \ge ... \ge \lambda_{m} \ge 0$ and $\lambda_{1},...,\lambda_{m} \in \mathbb{N} \cup \{ \infty \}$, then every orbit of the $G_{\infty}$ action on $\mathcal{M}_{\infty}$ has a unique following standard form with some $\lambda$:
	\begin{flalign*}
		\delta_{\lambda}=\begin{pmatrix}
			0      & ...      & 0      & t^{\lambda_{1}} & 0               & ...     & 0 \\
			0      & ...      & 0      & 0               & t^{\lambda_{2}} & ...     & 0 \\
			\vdots & \ddots   & \vdots & \vdots          & \vdots          & \ddots  & \vdots \\
			0      & ...      & 0      & 0               & 0               & ...     & t^{\lambda_{m}}
		\end{pmatrix}.
	\end{flalign*}
	Here we use the convention that $t^{\infty}=0$ and $\lambda_{i}$ can be $\infty$. The orbit with the standard form $\delta_{\lambda}$ is denoted by $\mathcal{C}_{\lambda}$.
\end{proposition}
Now we associated every orbit of the $G_{\infty}$ action on $\mathcal{M}_{\infty}$ with a $\lambda$, so the next step is to study how to characterize $\lambda$.
\begin{definition}
	$\lambda=(\lambda_{1},...,\lambda_{s},...)$(finitely many or infinitely many terms) is called a pre-partition if $\lambda_{1},...,\lambda_{s},... \in \mathbb{N} \cup \{ \infty \}$ and $\lambda_{1} \ge ... \ge \lambda_{s} \ge ...$. Define a relation $\triangleleft$ on two pre-partitions $\lambda=(\lambda_{1},\lambda_{2},...)$ and $\mu=(\mu_{1},\mu_{2},...)$ as the following:
	\begin{flalign*}
		\lambda \triangleleft \mu \iff \lambda_{i}+\lambda_{i+1}+... \le \mu_{i}+\mu_{i+1}+..., \  \forall i.
	\end{flalign*}
	The set of pre-partitions is a poset under the relation $\triangleleft$.
\end{definition}

We can also define a relation $\le$ on two orbits $\mathcal{C}, \tilde{\mathcal{C}}$ as the following:
\begin{flalign*}
	\mathcal{C} \le \tilde{\mathcal{C}} \iff \mathcal{C} \subset \bar{\tilde{\mathcal{C}}},
\end{flalign*}
where $\bar{\tilde{\mathcal{C}}}$ means the closure of $\tilde{\mathcal{C}}$. The set of orbits also forms a poset under this relation. The crucial connection between these two different relations is the following proposition.
\begin{proposition}
	(Theorem $4.7$ in \cite{12}) 
	For any two pre-partition $\lambda$ and $\mu$, with $\lambda=(\lambda_{1},...,\lambda_{m}),\mu=(\mu_{1},...,\mu_{m})$, we have $\lambda \triangleleft \mu \iff \mathcal{C}_{\mu} \le \mathcal{C}_{\lambda}$.
\end{proposition}
From this proposition we know studying the poset of orbits is equivalent to studying the poset of pre-partitions. Now in order to give the irreducible components of $\mathscr{L}^{m,n}_{r,k}$, we first define the contact locus, which connects $\mathscr{L}^{m,n}_{r,k}$ with arc space.
\begin{definition}
	The contact locus $Cont^{p}(I_{r}^{m,n})$ is defined to be all the arc $\alpha$ in $\mathcal{M}_{\infty}$ such that $order_{\alpha}(f)\ge p, \forall f \in I_{r}^{m,n}$. This also implies that $Cont^{p}(I_{r}^{m,n})=\pi_{\infty,p-1}^{-1}(\mathscr{L}^{m,n}_{r,p-1})$, where $\pi_{\infty,p-1}: \mathcal{M}_{\infty} \longrightarrow \mathcal{M}_{p-1}$ is the truncation map. 
\end{definition}

After studying the poset of pre-partitions, we can get the following proposition.
\begin{proposition}
	(Proposition $4.9$ in \cite{12}) 
	For every irreducible component $\mathcal{C}$ in $Cont^{p}(I_{r}^{m,n})$, there exists a unique $\lambda$ such that $\mathcal{C}=\bar{\mathcal{C}_{\lambda}}$, where $\lambda=(d \ d \ ... \ d \ e \ 0 \ ... \ 0)$, satisfying $p=(a+1)d+e$, $0 \le e < d$; either $e=0, 0 \le a \le r-1$ or $e > 0, 0 \le a < r-1$ and there are totally $a+m-r+1$ many $d$ in $\lambda$.
\end{proposition}

After calculating the number of $\lambda$ satisfying the condition above, we can get the number of irreducible components of $Cont^{p}(I_{r}^{m,n})$, and also the number of components of $\mathscr{L}^{m,n}_{r,k}$:
\begin{theorem}
	(Theorem $4.12$, Corollary $4.13$ in \cite{12}) 
	When $r=1, m$, $Cont^{k+1}(I_{r}^{m,n}) \subset \mathcal{M}_{\infty}$ and $\mathscr{L}^{m,n}_{r,k}$ is irreducible. When $1<r<m$, $Cont^{k+1}(I_{r}^{m,n})$ and $\mathscr{L}^{m,n}_{r,k}$ have $k+2-\lceil \frac{k+1}{r} \rceil$ many  irreducible components. 
\end{theorem}

Using the characterization for irreducible components in \cite{12}, we can show that the principal component always exists, as described in the following proposition:
\begin{proposition}
In Proposition $8.5$, we take $d=p,a=0,e=0$, when $\lambda$ is a combination of $m-r+1$ many $d$ and $r-1$ many $0$. Considering the projection map $\pi_{\infty}:\mathcal{M}_{\infty} \longrightarrow \mathcal{M}$, then we have $\bar{\mathcal{C}_{\lambda}}=\bar{\pi_{\infty}^{-1}((\mathscr{L}^{m,n}_{r})_{reg})} \cap Cont^{p}(I_{r}^{m,n})$, where $(\mathscr{L}^{m,n}_{r})_{reg}$ denotes the smooth part of $\mathscr{L}^{m,n}_{r}$.	
\end{proposition}
\begin{proof}
For traditional determinantal varieties, we have $(\mathscr{L}^{m,n}_{r})_{sing}=\mathscr{L}^{m,n}_{r-1}$, i.e., the singularity of a determinantal variety defined by an $r$-th minors is exactly the part of it that is defined by an $r-1$-th minors. For the point $P$ in $\mathcal{C}_{\lambda}$, consider the ideal generated by the $r-1$-th minors, which is the same as the ideal generated by the $r-1$-th minors of $\delta_{\lambda}$, since the group action does not change it. This tells us that there exists an $r-1$-th minor invertible, i.e., $P \in \pi_{\infty}^{-1}((\mathscr{L}^{m,n}_{r})_{reg})$, and hence $\bar{\mathcal{C}_{\lambda}} \subset \bar{\pi_{\infty} ^{-1}((\mathscr{L}^{m,n}_{r})_{reg})}$. Combined with $\mathcal{C}_{\lambda} \subset Cont^{p}(I^{m,n}_{r})$ it is known that $\bar{\mathcal{C}_{\lambda}} \subset \bar{\pi_{\infty}^{-1}((\mathscr{L}^{m,n}_{r})_{reg})}$. Conversely, for the point $P$ in $\pi_{\infty}^{-1}((\mathscr{L}^{m,n}_{r})_{reg})\cap Cont^{p}(I_{r}^{m,n})$, let $P \in \mathcal{C}_{\mu}$, where $\mu=(\mu_{1}, ... ... ,\mu_{m})$ is some pre-partition. Then, from the fact that the constant term of $P$ is not in the ideal generated by the $r-1$-th minors, we know that $\mu_{m}=... =\mu_{m-r+2}=0$, and by $P \in Cont^{p}(I_{r}^{m,n})$ we know that $\mu_{m-r+1} \ge p$, which tells us that $\lambda \triangleleft \mu$, and thus $P \in \mathcal{C}_{\mu} \subset \bar{\mathcal{C}_{\lambda}}$, and the above proposition is proved.
\end{proof}

The above proposition in fact describes the jet schemes of determinantal varieties by lifting it to arc space, and by composing with $\pi_{\infty,p-1}$ on both sides of the conclusion of the above proposition we get: $\pi_{\infty,p-1}(\bar{\mathcal{C}_{\lambda}})=\bar{\pi_{p-1 }^{-1}((\mathscr{L}^{m,n}_{r})_{reg})}$, i.e., the closure of inverse image of the smooth part of the bottom space. This tells us that $\pi_{\infty,p-1}(\bar{\mathcal{C}_{\lambda}})$ is in fact the principal component of $\mathscr{L}^{m,n}_{r,p-1}$. What was done in \cite{5} is to work out the Hilbert series of the principal component of $\mathscr{L}^{m,n}_{2,1}$, and to show that it is the square of the Hilbert series of $\mathscr{L}^{m,n}_{2}$; and what we did in Sections $4,5,6$ and $7$ is to work out the Hilbert series of the $\mathscr{L }^{3,n}_{3,1}$ comes to a similar conclusion, so for the general $\mathscr{L}^{m,n}_{r,k}$ we have the following conjecture.
\begin{conjecture}
Denote $\tilde{I^{m,n}_{r,k}}$ as the ideal corresponding to the principal component of $\mathscr{L}^{m,n}_{r,k}$, then its corresponding abstract simplicial complex is shellable, and the Hilbert series of the principal component of $\mathscr{L}^{m,n}_{r,k}$ is exactly $k+1$-th square of the Hilbert series of $\mathscr{L}^{m,n}_{r}$.
\end{conjecture}

The characterization for irreducible components in \cite{12} is reduced to the closure of the orbits, which is a geometric characterization but does not directly give the ideal of an irreducible component, and we next need to find a way to calculate the specific polynomials that define each irreducible component. Our approach is to compute minors of different sizes and try to distinguish between different orbits by these minors in order to draw our conclusions. First we give the following notation: for the arc space $\mathcal{M}_{\infty}$, let the elements of the $i$-th row and $j$-th column be $x_{i,j}^{(0)}+x_{i,j}^{(1)}t+... $, and we denote that the family of polynomials consisting of the $t^{\alpha}$ terms of its $s \times s$ minors is $\Theta_{s,\alpha}$. We have the following proposition.

\begin{proposition}
Let $1<r<m$, and for each irreducible component $\bar{\mathcal{C}_{\lambda}}$ in $Cont^{p}(I^{m,n}_{r})$, we define the following family of polynomials: for $1 \le s \le m$, suppose that the ideal generated by the $s$-th minors of $\delta_{\lambda}$ is $(t^{i_{s}})$ (in fact $i_{s}=\lambda_{m}+... +\lambda_{m-s+1}$), then we take $\Theta_{s,0},... ,\Theta_{s,i_{s}-1}$ into this family of polynomials. We do the above procedure for each $s$ to get the whole family of polynomials and denote it as $\Omega_{\lambda}$. Then the set of zeros defined by $\Omega_{\lambda}$ is exactly $\bar{\mathcal{C}_{\lambda}}$.
\end{proposition}
\begin{proof}
For any point $P$ in $\mathcal{C}_{\lambda}$, the ideal of $P$ generated by the $s \times s$ minors is the same as the ideal generated by the $s \times s$ minors of $\delta_{\lambda}$ because the group action does not change it. This tells us that $P$ satisfies the polynomials in $\Omega_{\lambda}$, so we have $\bar{\mathcal{C}_{\lambda}} \subset Z(\Omega_{\lambda})$, where $Z(\Omega_{\lambda})$ denotes the $\Omega_{\ lambda}$ the set of zeros.

Conversely, for each point $P$ that satisfies polynomials in $\Omega_{\lambda}$, we assume that $P \in \mathcal{C}_{\mu}$, where $\mu=(\mu_{1}, ... ,\mu_{m})$ is a pre-partition. Then the standard type of $P$ is $\delta_{\mu}$. Since $P$ satisfies the equations in $\Omega_{\lambda}$, we have $\mu_{m}+\mu_{m-1}+... +\mu_{m-s+1} \ge i_{s}$ holds for all $1\le s \le m$. This tells us $\lambda \triangleleft \mu$, and by Proposition $8.3$, this is equivalent to $\mathcal{C}_{\mu} \subset \bar{\mathcal{C}_{\lambda}}$, so we have $P \in \mathcal{C}_{\mu} \subset \bar{\mathcal{C}_{\lambda}}$.
\end{proof}

The $\Omega_{\lambda}$ we get contains a very large number of polynomials, but in fact we don't need that many; the polynomials in a part of $\Omega_{\lambda}$ are enough to get the whole of $\Omega_{\lambda}$, as described in the following proposition.
\begin{proposition}
For each $1 \le s \le m-1$, let $j_{s}=\lceil\frac{i_{s}}{s}\rceil$, and we define a new family of polynomials, $\Omega_{\lambda}'$, as follows: first we take $\Theta_{1,0}, ... ,\Theta_{1,i_{1}-1} \in \Omega_{\lambda}'$; when $i_{s}+j_{s}<i_{s+1}$, take $\Theta_{s+1,i_{s}+j_{s}},... ,\Theta_{s+1,i_{s+1}-1} \in \Omega_{\lambda}'$. This gives us the subfamily $\Omega_{\lambda}'$ of $\Omega_{\lambda}$, and then $\Omega_{\lambda}'$ can generate $\Omega_{\lambda}$.
\end{proposition}

\begin{proof}
Let $\Omega_{\lambda,q}=\{\Theta_{j,l} \in \Omega_{\lambda}|j \le q\}$ and $\Omega_{\lambda,q}'=\{\Theta_{j,l} \in \Omega_{\lambda}'|j \le q\}$. We prove inductively for $q$ that $\Omega_{\lambda,q}$ can be generated by $\Omega_{\lambda,q}'$. When $q=1$, $\Omega_{\lambda,1}'=\Omega_{\lambda,1}$ and the assertion holds. Now assume that the assertion holds for $q<q_{0}$ and consider the case $q=q_{0}$. $\Omega_{\lambda,q_{0}}$ has more $\Theta_{q_{0},0},...$ compared to $\Omega_{\lambda,q_{0}-1}$. ,$\Theta_{q_{0},i_{q_{0}}-1}$.
	
When $i_{q_{0}-1}+j_{q_{0}-1}<i_{q_{0}}$, it is sufficient to prove that $\Theta_{q_{0},0},... ,\Theta_{q_{0},i_{q_{0}-1}+j_{q_{0}-1}-1}$ can be generated by $\Omega_{\lambda,q_{0}-1}$. Indeed, for any $f \in \Theta_{q_{0},l}(l \le i_{q_{0}-1}+j_{q_{0}-1}-1)$, $f$ is the $l$-secondary term in the $q_{0}$-order determinant expansion without loss of generality, and we assume that the determinant in which $f$ is located is in the first $q_{0}$ rows and the first $q_{0}$ columns. Then $f$ has the following expansion:
\begin{flalign*}
	f=\sum_{l_{1}+...+l_{q_{0}}=l}\begin{vmatrix}
		x_{1,1}^{(l_{1})} & x_{1,2}^{(l_{1})} & \cdots  & x_{1,q_{0}}^{(l_{1})} \\
		x_{2,1}^{(l_{2})} & x_{2,2}^{(l_{2})} & \cdots  & x_{2,q_{0}}^{(l_{2})} \\
		\vdots            &   \vdots          & \ddots  & \vdots                \\
		x_{q_{0},1}^{(l_{q_{0}})} & x_{q_{0},2}^{(l_{q_{0}})} & \cdots & x_{q_{0},q_{0}}^{(l_{q_{0}})}
	\end{vmatrix}.
\end{flalign*}	
For each term in its expansion as above, we take $l_{i}=max\{l_{1}, ... ,l_{q_{0}}\}$ and expand it by the $i$-th row, so that it can be obtained after computation:
\begin{flalign*}
	f=\sum_{l_{i}=max\{l_{1},...,l_{q_{0}}\},i}\sum_{l_{i}}\sum_{l_{1}+...+l_{q_{0}}-l_{i}=l-l_{i}}\begin{vmatrix}
		x_{1,1}^{(l_{1})} & x_{1,2}^{(l_{1})} & \cdots  & x_{1,q_{0}}^{(l_{1})} \\
		x_{2,1}^{(l_{2})} & x_{2,2}^{(l_{2})} & \cdots  & x_{2,q_{0}}^{(l_{2})} \\
		\vdots            &   \vdots          & \ddots  & \vdots                \\
		x_{q_{0},1}^{(l_{q_{0}})} & x_{q_{0},2}^{(l_{q_{0}})} & \cdots & x_{q_{0},q_{0}}^{(l_{q_{0}})}
	\end{vmatrix}.
\end{flalign*}
When fixing $i,l_{i}$, what is summed under the last summation number can be generated by elements in $\Theta_{q_{0}-1,l-l_{i}}$ after expanding by the $i$-th row. We assert that $l-l_{i} \le i_{q_{0}-1}-1$. In fact, if there is $l-l_{i} \ge i_{q_{0}-1}$, since $l_{i} \ge \lceil\frac{l}{q_{0}}\rceil$, $l-l_{i} \le l-\lceil\frac{l}{q_{0}}\rceil \le i_{q_{0}-1}+j_{q _{0}-1}-1-\lceil\frac{i_{q_{0}-1}+j_{q_{0}-1}-1}{q_{0}}\rceil$, so $j_{q_{0}-1} \ge 1+ \lceil\frac{i_{q_{0}-1}+j_{q_{0}-1}-1}{q_{0}}\rceil \ge 1+\frac{i_{q_{0}-1}+j_{q_{0}-1}-1}{q_{0}}$. Calculation gives $(q_{0}-1)\lceil\frac{i_{q_{0}-1}}{q_{0}-1}\rceil \ge i_{q_{0}-1}+q_{0}-1$, a contradiction. Thus $l-l_{i} \le i_{q_{0}-1}-1$, which tells us that $f$ can be generated by $\Omega_{\lambda,q_{0}-1}$.

When $i_{q_{0}-1}+j_{q_{0}-1} \ge i_{q_{0}}$, it can be similarly shown that $\Theta_{q_{0},l}$ can be generated by $\Omega_{\lambda,q_{0}-1}$, where $l \le i_{q_{0}-1}+j_{q_{0}-1}-1$, and so $l \le i_{q_{0}}-1$, and thus $\Omega_{\lambda,q}'$ can be generated by $\Omega_{\lambda,q}$, and the proof is complete.
\end{proof}
Notice that the projection mapping of $\pi_{\infty,p-1}:Cont^{p}(I^{m,n}_{r}) \longrightarrow \mathscr{L}^{m,n}_{r,p-1}$ actually replaces all of the variants $x_{i,j}^{(l)}(l \ge p)$ by $0$, and so for $\Theta_{s,\alpha}$, we will denote the corresponding family of polynomials obtained by replacing all the variables $x_{i,j}^{(l)}(l \ge p)$ by $0$, so for the elements of $\Theta_{s,\alpha}$ we will denote the corresponding family of polynomials obtained by constructing $\Theta_{s,\alpha}$ in the same way as in Proposition $8.10$ as $\tilde{\Theta}_{s,\alpha}$. so that we obtain the family of polynomials defining the irreducible components of $\mathscr{L}^{m,n}_{r,p-1}$, i.e., Theorem \ref{mthc}.
\begin{theorem}
$\tilde{\Omega_{\lambda}'}$ can be used as a family of polynomials defining the irreducible component $\pi_{\infty,p-1}(\bar{\mathcal{C}_{\lambda}})$.
\end{theorem}

\section{Acknowledgement}
On behalf of all authors, the corresponding author states that there is no conflict of interest.

\end{document}